\documentclass[twocolumn,aps,floatfix,superscriptaddress,longbibliography]{revtex4-2}
\pdfoutput=1 
\usepackage{amsmath,amssymb,eucal,graphicx,float,epstopdf,xparse}
\usepackage{booktabs}
\usepackage{epsfig,subfigure}
\usepackage[utf8]{inputenc}

\NewDocumentCommand{\eulerian}{omm}
 {%
  \genfrac<>{0pt}{}{#2}{#3}%
  \IfValueT{#1}{_{\!#1}}%
 }

\newcommand{\rev}[1]{{#1}}
\usepackage{xcolor}
\newcommand{\revv}[1]{#1}

\usepackage[colorlinks=true, urlcolor=blue, anchorcolor=blue, citecolor=blue,filecolor=blue,linkcolor=blue,menucolor=blue]{hyperref}
\usepackage[capitalize]{cleveref}

\begin{document}
\title{Random Recursive Simplicial Complexes}

\author{P. L. Krapivsky}
\affiliation{Department of Physics, Boston University, Boston, Massachusetts 02215, USA}
\affiliation{Santa Fe Institute, Santa Fe, New Mexico 87501, USA}

\author{M. Lucas}
\affiliation{Department of Mathematics and Namur Institute for Complex Systems (naXys), Université de Namur, Namur, Belgium}
\affiliation{Mycology Laboratory, Earth and Life Institute, Université Catholique de Louvain, Louvain-la-Neuve, Belgium}

\begin{abstract} 
We investigate random recursive simplicial complexes that grow by adding, at each step, a vertex together with a simplex formed by joining the new vertex to a randomly chosen existing simplex. We also add all faces of the new simplex to ensure that the resulting object remains a simplicial complex. If the choice of an existing simplex is uniform among simplices of dimension $<m$, the total number $S_d$ of simplices of any admissible dimension $d\leq m$ is an asymptotically self-averaging random variable. \revv{Mapping the growth onto a generalized P\'olya urn, we prove this for every $m$ and determine the asymptotic growth law of the average of $S_d$ from a linear eigenvalue problem. The fluctuations of $S_d$ are Gaussian for $m\le5$ and anomalously large for $m\ge6$.} We also probe the degree distribution, \revv{whose algebraic tail we characterize for every $m$,} examine the probabilities of various extreme outcomes, and analyze the characteristics of the first vertex. 
\end{abstract}

\maketitle

\section{Introduction}
\label{sec:intro}

\revv{Simplicial complexes and related objects, such as CW complexes constructed by gluing balls together, play crucial roles in algebraic topology \cite{Hatcher} and combinatorics \cite{Flajolet}.} A 0-dimensional simplex is a point; a 1-dimensional simplex is a segment; a 2-dimensional simplex is a triangle; a 3-dimensional simplex is a tetrahedron; etc. A $(d-1)$-dimensional simplex has $2^d-1$ non-empty subsimplices, faces in short: $\binom{d}{1}$ vertices, $\binom{d}{2}$ edges, $\binom{d}{3}$ triangles, etc. A simplicial complex $\mathcal{K}$ is a set of simplices such that every face of a simplex from $\mathcal{K}$ is also in $\mathcal{K}$ and the non-empty intersection of any two simplices from $\mathcal{K}$ is a face of both simplices \cite{Hatcher}. The dimension of a simplicial complex is the maximal dimension of its simplices. One-dimensional simplicial complexes are often referred to as graphs. In graph theory \cite{Diestel}, points are called vertices and segments are called edges. 

\rev{\revv{Beyond their role in mathematics, simplicial complexes are increasingly used as natural models of complex systems} with interactions between more than two units~\cite{battiston2020networks,bick2023what}---for example in neural, social, and collaboration networks---which can reshape collective dynamics \cite{iacopini2019simplicial,skardal2019abrupt,zhang2023higherorder,battiston2026collective,millan2025topology}. This has motivated a search for tractable generative models of higher-order structure. Studies of large random simplicial complexes are rather recent; see, e.g., \cite{Pippenger,3-manifolds,Linial06,Meshulam,Linial16, Farber16,Pittel16,Ginestra16-NGF,Farber17a,Farber17b,Bobrowski17, Ginestra17,Ginestra-Sergey,Ginestra18-NGF,Ginestra-Sergey20,Kahle,Petri-B21, Ginestra21,Petri-B22,Dima22}. 
Among these, growth models fall into two categories. Deterministic ones build the complex by a recursive rule---such as the pseudofractal simplicial complexes \cite{Sergey02}---yielding scale-free structures with exactly computable properties \cite{nurisso2025higherorder}. Stochastic ones add simplices at random: examples include the \revv{network geometry with flavor} \cite{Ginestra16-NGF} and its dimension-wise variant \cite{febbe2026model}, which grow $d$-dimensional complexes through an effective preferential attachment. The model introduced below is stochastic but, unlike these, grows on a uniformly chosen simplex---the higher-dimensional counterpart of the recursive random tree.
}

In this paper, we investigate random recursively growing simplicial complexes. We begin the process with a single `primordial' vertex: $\mathcal{K}_1=\langle\{v_1\}\rangle$. On the $n$th step, we add a new vertex $\{v_n\}$ and join it with randomly chosen simplex from $\mathcal{K}_{n-1}$, say $\{v_{i_1},\ldots,v_{i_k}\}\in \mathcal{K}_{n-1}$ to form the simplex $\{v_{i_1},\ldots,v_{i_k},v_n\}$. We add this simplex together with all its faces to $\mathcal{K}_{n-1}$ to form $\mathcal{K}_n$.

\begin{figure}[b]
	\centering
	\includegraphics[width=1\linewidth]{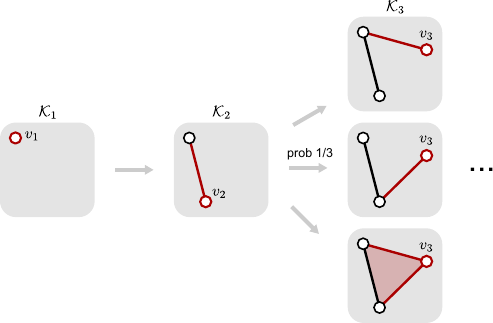}
	\caption{Illustration of the RRSC(2) model: the first three steps yield simplicial complexes $\mathcal{K}_1$, $\mathcal{K}_2$, and $\mathcal{K}_3$. At each step, a new vertex and simplices (red) are added at random. Each of the three possible $\mathcal{K}_3$ occurs with probability $1/3$.}
	\label{fig:growth}
\end{figure}

The above null model generates growing random simplicial complexes. The dimension of a typical simplicial complex with $N$ vertices diverges as $N\to\infty$. Simplicial complexes of a fixed dimension usually arise in applications \cite{Diestel,Hatcher,Flajolet}. A natural deformation of the above rules allows the generation of simplicial complexes of any fixed dimension. Suppose a simplex joining an arriving vertex is chosen uniformly among $0$-dimensional simplices, i.e., among vertices. The emerging random simplicial complexes are one-dimensional graphs, more precisely, trees, as there is exactly one path connecting any two vertices. This recursive random tree (RRT) model is a paradigmatic parameter-free model of growing random trees, see, e.g., \cite{Pittel94,KR01,KR02-fluct,KR02,Janson05,Drmota,Hofstad,Janson15,Janson19}. Similarly, suppose an arriving vertex joins with (uniformly chosen) existing simplex of dimension smaller than $m$. This rule defines the RRSC(m) model, where `RRSC' stands for `random recursive simplicial complex'. Thus, we have a class of models parametrized by positive integers $m \in \mathbb{Z}_+$. The RRSC(1) model is identical to the RRT model. \rev{The first three steps of the RRSC(2) model are illustrated in \cref{fig:growth}.} \revv{The RRSC model resembles the random recursive hypergraph (RRH) model \cite{PK-hyper}. A class of RRH(m) models, which are hypergraph analogs of the RRSC(m) models, is also easy to define.}

\revv{Systems with higher-order interactions, such as collaborations, involve groups of many sizes~\cite{battiston2020networks,Kleinberg22b}. A newcomer may join a single collaborator, a pair, or an existing team. In a simplicial complex, every subgroup of a group is itself a face~\cite{Hatcher}, so each of these faces is an equally valid entry point. The RRSC($m$) rule chooses among them uniformly, with $m+1$ as the largest group size and no other parameter. Growth by attachment to top faces only~\cite{Ginestra16-NGF,FIMS22,DMS01} (our HSC($m$) model, Appendix~\ref{ap:hom}) instead lets newcomers join full-size groups only. All maximal simplices then have dimension $m$, and the numbers of simplices are deterministic. Mixing dimensions makes these numbers random and gives the degree of a vertex a memory of the faces that contain it, which changes the degree exponent. The RRSC($m$) models are thus the simplest growth models of simplicial complexes with mixed dimensions, and, as we show below, they remain tractable: the numbers of simplices can be computed rigorously for every $m$, while the degree statistics display features absent from growth by attachment to top faces.}

In Sec.~\ref{sec:basic}, we compute the probability $T_N$ that the emerging simplicial complex is one-dimensional, i.e., a tree in the terminology of graph theory \cite{Diestel}. We show that 
\begin{equation}
\label{TN-sol}
T_N = \frac{(N-1)!}{(2N-3)!!}
\end{equation}
for all $m\geq 2$. We also compute a few similar probabilities of extreme outcomes, such as the emergence of a star or a linear graph. The results of Sec.~\ref{sec:basic} are exact, i.e., they hold for an arbitrary total number of vertices $N$. 

In Sec.~\ref{sec:simplices}, we analyze the total number of simplices $S_d$ of every feasible dimension ($d=0,1,2,\ldots, m$) in simplicial complexes generated by the RRSC(m) model. In an ensemble of simplicial complexes generated by the RRSC(m) model, $S_d$ are random quantities that become asymptotically self-averaging when $N\gg 1$, so the average values provide the chief information. We also \revv{show} extensivity, i.e., a linear growth in $ N$ of the average values $\big\langle S_d \big\rangle$\revv{, rigorously via the urn mapping of Sec.~\ref{sec:urn}}. For the RRSC(m) model, only $S_d$ with $d=0,1,2,\ldots, m$ are non-trivial. (The more precise $S_d(m,N)$ notation is often shorten to  $S_d$.)

The quantities $S_d$ are random except for the total number of vertices
\begin{subequations}
\begin{equation}
\label{S0}
S_0 = S_0(m,N) = N
\end{equation}
The alternating sum of $S_d$, the Euler characteristic, is also deterministic and equals to one 
\begin{equation}
\label{Euler:def}
\chi = \sum_{d=0}^m (-1)^d S_d = 1
\end{equation}
\end{subequations}
Indeed, the definition of the RRSC(m) model implies that the emerging simplicial complexes are topologically trivial \cite{Hatcher}, viz., contractible to a point. This observation explains \eqref{Euler:def}. 

Relations \eqref{S0}--\eqref{Euler:def} show that for the RRSC(1)=RRT model, both $S_0$ and $S_1=N-1$ are deterministic. For RRSC(2) model, only one of the two random variables $S_1$ and $S_2$ is independent. On average, the variables $S_d$ exhibit an extensive growth:
\begin{equation}
\label{Sd:av}
\big\langle S_d \big\rangle \simeq a_d(m)N
\end{equation}
when $N\gg 1$. The variances $\text{Var}[S_d]=\langle S^2_d \rangle-\langle S_d \rangle^2$ also grow linearly in the system size \revv{for $m\le5$}:
\begin{equation}
\label{Sd:var}
\text{Var}[S_d] \simeq v_d(m)N
\end{equation}
when $N\gg 1$. \revv{The variables $S_d$ are asymptotically Gaussian when $m\le5$; for $m\ge6$ their fluctuations grow faster than $\sqrt{N}$ (Sec.~\ref{sec:urn}). These statements follow from a mapping of the RRSC($m$) process onto a generalized P\'olya urn, under which the amplitudes $a_d(m)$ become the components of the Perron eigenvector of an $m\times m$ integer matrix.} \revv{We give explicit derivations for small $m$; the general case follows from the urn mapping of Sec.~\ref{sec:urn}.} Specifically, we analyze the RRSC(2) model in detail in Appendix \ref{ap:edges}. We show that the number of edges $S_1$ is an asymptotically Gaussian random variable and determine the asymptotic behavior of the variance. We then outline a generalization to the RRSC(3) model, compute the covariance matrix of $S_1$ and $S_2$, and show that the probability distribution $P(S_1,S_2; N)$ is asymptotically Gaussian. In Sec.~\ref{sec:simplices}, we establish lower and upper bounds for $S_d$ that are asymptotically linear in $N$, thereby corroborating the extensivity of $S_d$. 

Since the variables $S_d$ become self-averaging in the $N\to\infty$ limit, we are mostly interested in the growth of the averages, i.e., in computing the amplitudes $a_d(m)$ in \eqref{Sd:av}. In Sec.~\ref{sec:linear}, we derive a system of algebraic equations for $a_d(m)$ and establish explicit expressions for the amplitudes in the case of $m=2$. When $m\geq 3$, one can numerically solve the algebraic equations and obtain the solutions with any desired accuracy. \revv{We show that these equations reduce to a linear eigenvalue problem.} In Appendix \ref{ap:10}, we present the amplitudes for the RRSC(m) models with $m\leq 9$. 

In Sec.~\ref{sec:primordial}, we consider the primordial vertex and discuss the characteristics of the RRSC(m) models associated with the primordial vertex. We show that the degree of the primordial vertex scales algebraically with $N$ and is not self-averaging\revv{; the urn mapping explains this through a random prefactor set by the early growth of the neighbourhood of the vertex}. The average degree of the primordial vertex should be comparable to the maximal degree. Relying on this observation and assuming that the degree distribution has an algebraic tail, $n_\delta\sim \delta^{-\nu}$ for $\delta\gg 1$, we determine the decay exponent $\nu(m)$ \revv{for $m\le4$ (e.g., $\nu(3)\approx2.715\,357$, $\nu(4)\approx2.449\,639$). Using the urn mapping, we then show, without this assumption, that for general $m$, $\nu(m)=1+\lambda_1(m)/\lambda_1(m-1)$, where $\lambda_1(m)$ is the largest eigenvalue of an $m\times m$ integer matrix.} For the RRSC(2) model, the decay exponent has a simple explicit form: 
\begin{equation}
\label{tail-2}
\nu(2) = \frac{5+\sqrt{5}}{2} = 3.618\,033\,988\ldots
\end{equation}

In Sec.~\ref{sec:degree}, we outline the scheme that may allow the determination of the degree distribution. The procedure is recursive in character. The amount of computations rapidly increases with degree, and we only determined the fraction of vertices of degree one, which we expressed via the amplitudes $a_1(m),\ldots,a_{m-1}(m)$. The prediction has a simple, explicit form for the RRSC(2) model
\begin{equation}
\label{deg:1-2}
n_1 = \frac{7-\sqrt{5}}{22}
\end{equation}

\revv{Our results differ in level of rigor, as summarized in Table~\ref{tab:status}. Those of Secs.~\ref{sec:basic} and~\ref{sec:simplices}, the probabilities of extreme outcomes and the bounds on $S_d$, are exact for all $N$. The equations for the amplitudes in Sec.~\ref{sec:linear} are derived assuming that the $S_d$ are self-averaging. The urn mapping of Sec.~\ref{sec:urn} proves this assumption, so the amplitudes and the nature of the fluctuations of $S_d$ are rigorous as $N\to\infty$. The variances (Appendix~\ref{ap:edges}), the growth of the average degree of the primordial vertex, and the fraction $n_1$ of vertices of degree one rely on mean-field approximations, and the exponent $\nu$ on heuristic arguments. We support them with numerical simulations.}

The degree of a vertex counts the number of adjacent edges and provides a natural local characteristic for graphs. The degree of a vertex that we used in Sec.~\ref{sec:degree} is defined as the degree of the one-dimensional graph skeleton of the simplicial complex. One can extend the concept of degree from graphs to simplicial complexes by counting $D$-dimensional simplices adjacent to a $ d$-dimensional simplex. Such generalized degrees are parametrized by pairs $(d, D)$ with $0\leq d<D$; the $(0,1)$ degree is the standard one. For the RRSC(m), there are $m(m+1)/2$ degree distributions. 

In Sec.~\ref{sec:future}, we discuss directions for future work. In particular, we define another natural class of HSC(m) models generating random {\em homogeneous} simplicial complexes with all maximal simplices having the same dimension $m$. In Appendix~\ref{ap:hom}, we show that the HSC(m) models are more tractable than the RRSC(m) models.

\section{Extremal characteristics}
\label{sec:basic}

The primordial complex is $\mathcal{K}_1=\langle\{v_1\}\rangle$. On the next step, $N=2$, the complex is still deterministic:
\begin{equation}
\label{K-2}
\mathcal{K}_2 = \langle\{v_1\}, \{v_2\}; \{v_1, v_2\}\rangle
\end{equation}
An alternative shorter notation $[\{v_1, v_2\}]$ for the simplicial complex \eqref{K-2} is equally informative since all faces of each simplex belong to the simplicial complex, and we can recover the detailed description. 

Starting from $N=3$, more than one simplicial complex can be built. For $N=3$, three simplicial complexes occur with equal probabilities:
\begin{equation}
\label{K3}
\begin{cases}
 [\{v_1, v_2\}, \{v_1, v_3\}] & \text{prob} ~\frac{1}{3} \\
 [\{v_1, v_2\}, \{v_2, v_3\}] & \text{prob} ~\frac{1}{3}\\
 [\{v_1, v_2,v_3\}] & \text{prob} ~\frac{1}{3}
\end{cases}
\end{equation}
for the RRSC(m) models with $m\geq 2$, \rev{as illustrated in \cref{fig:growth}}. In the detailed notation, the outcomes \eqref{K3} are
\begin{equation*}
\begin{split}
&\langle\{v_1\}, \{v_2\},  \{v_3\}; \{v_1, v_2\}, \{v_1, v_3\}\rangle \\
&\langle\{v_1\}, \{v_2\}, \{v_3\}; \{v_1, v_2\}, \{v_2, v_3\}\rangle \\
&\langle\{v_1\}, \{v_2\}, \{v_3\}; \{v_1, v_2\}, \{v_1, v_3\}, \{v_2, v_3\}; \{v_1, v_2,v_3\}\rangle
\end{split}
\end{equation*}

By definition, the dimensions of simplicial complexes built via the RRSC(m) procedure lie in the range
\begin{equation}
\label{dim}
1\leq \text{dim}(\mathcal{K})\leq m \qquad\text{when}\quad N\geq 2
\end{equation}

The primordial simplicial complex $\mathcal{K}_1=\langle\{v_1\}\rangle$ is zero-dimensional. The next simplicial complex \eqref{K-2} is one-dimensional: $\text{dim}(\mathcal{K}_2)=1$. Simplicial complexes of the maximal dimension may be built when $N\geq m+1$. One anticipates that for $N\gg m$, the typical simplicial complex built via the RRSC(m) procedure has dimension $m$. What is the probability of building a simplicial complex of the minimal dimension, $\text{dim}(\mathcal{K})=1$? Such a simplicial complex is a tree with $N$ vertices and $N-1$ edges. Denote by $T_N(m)$ the probability of this outcome:
\begin{equation}
\label{TN-def}
T_N(m) = \text{Prob}[\text{dim}(\mathcal{K})=1]
\end{equation}
For the RRSC(1), the simplicial complexes are trees by construction: $T_N(1)=1$. For the RRSC(m) models with $m\geq 2$, the probabilities $T_N$ are non-trivial but independent of $m$.  The recursive nature of the process implies 
\begin{equation}
T_{N+1} = \frac{N}{2N-1}\,T_N
\end{equation}
for $N\geq 2$. Solving this recurrence subject to the `initial condition' $T_2=1$ yields the announced probability \eqref{TN-sol}. This probability decreases with $N$ as $T_N\simeq 2^{-N}\sqrt{4\pi N}$.

There are many simplicial complexes of size $N$ with minimal dimension. Two special trees are stars and linear graphs. The star centered at the primordial vertex 
\begin{equation}
\label{star}
[\{v_1,v_2\}, \{v_1,v_3\},\ldots, \{v_1,v_N\}] 
\end{equation}
is formed with probability $\sigma_N$ found from the recurrence
\begin{equation}
\label{sigma:rec}
\sigma_{N+1} = \frac{1}{2N-1}\,\sigma_N
\end{equation}
Iterating \eqref{sigma:rec} starting from $\sigma_2=1$ yields
\begin{equation}
\label{sigma-sol}
\sigma_{N+1} = \frac{1}{(2N-1)!!}=\frac{2^N\,N!}{(2N)!}
\end{equation}
The second star
\begin{equation}
\label{star-2}
[\{v_1,v_2\}, \{v_2,v_3\}, \{v_2,v_4\}, \ldots, \{v_2,v_N\}] 
\end{equation}
is formed with the same probability $\sigma_N$. 

The linear graph
\begin{subequations}
\begin{align}
[\{v_1,v_2\}, \{v_2,v_3\}, \{v_3,v_4\}, \ldots, \{v_{N-1},v_N\}] 
\end{align}
with root at $v_1$ is formed with probability $\sigma_N$. The dual linear graph 
\begin{equation}
[\{v_1,v_2\}, \{v_1,v_3\}, \{v_3,v_4\}, \ldots, \{v_{N-1},v_N\}] 
\end{equation}
\end{subequations}
with root at $v_2$ is also formed with probability $\sigma_N$. The probability $L_N$ to form an arbitrary linear graph is found from the recurrence
\begin{equation}
\label{LN:rec}
L_{N+1} = \frac{2}{2N-1}\,L_N
\end{equation}
subject to the `initial' condition $L_2=1$. Thus 
\begin{equation}
\label{LN-sol}
L_{N} = \frac{2^{N-2}}{(2N-3)!!}
\end{equation}
All above probabilities are independent of $m\geq 2$, so we used the shorthand notations $\sigma_N$ and $L_{N}$.

\section{Simplices of various dimensions}
\label{sec:simplices}

A simplicial complex $\mathcal{K}$ with $N$ vertices built via the RRSC(m) procedure admits a decomposition
\begin{equation}
\label{decomp}
\mathcal{K} = \left\langle \mathcal{K}_0; \ldots; \mathcal{K}_m\right\rangle
\end{equation}
where $\mathcal{K}_d$ is the set of simplices in $\mathcal{K}$ of dimension $d$. All vertices $\{v_1\}, \ldots, \{v_N\}$ appear in $\mathcal{K}_0$. We shortly write $S_d$ for the number of simplices of dimension $d$ in $\mathcal{K}$. The quantities $S_d = |\mathcal{K}_d|$ depend on $\mathcal{K}$, the dimension $m$ and $N$. We shortly denote these random quantities by $S_d$. 

The total number of simplices 
\begin{equation}
\label{SN:def}
S = \sum_{d=0}^m S_d
\end{equation}
depends on $\mathcal{K}$. Therefore, $S$ is a random quantity (which also depends on $m$ and $N$). 

The number of vertices, $S_0=N$, is deterministic and independent of $m$. The alternating sum of the random quantities $S_d$ is the Euler characteristic \cite{Hatcher}. With our definition of the RRSC(m) procedure and the initial condition, the Euler characteristic is equal to one independently of $m$ and $N$, as we asserted in \eqref{Euler:def}. 

For the RRSC(1) model generating trees, all quantities $S_0, S_1,S$ are deterministic:
\begin{equation}
\label{S01}
S_0 = N, \quad S_1 = N-1, \quad S=2N-1
\end{equation}

\begin{figure}[t]
	\begin{center}
		\includegraphics[width=0.5\linewidth]{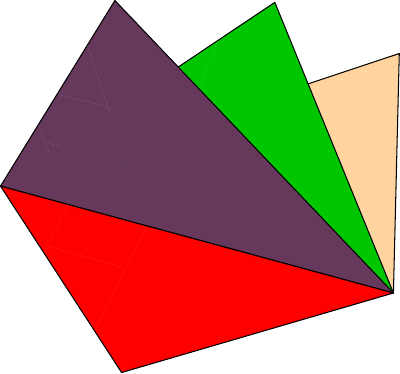}
		\caption{An open book simplicial complex with four triangular pages.} 
		\label{fig:book_illust}
	\end{center}
\end{figure}

For the RRSC(2), the quantities $S_1,S_2,S$ are random:
\begin{equation}
\label{SN12:bounds}
\begin{split}
N-1 &\leq S_1 \leq 2N-3  \\
  0  &\leq S_2 \leq N-2  \\
2N-1 & \leq S \leq 4N-5
\end{split}
\end{equation}
The lower bound is achieved on trees. The upper bound is achieved on simplicial complexes with maximal number of triangles such as an open book complex 
\begin{align}
\label{book}
 [\{v_1, v_2,v_3\}, \{v_1, v_2,v_4\}, \ldots, \{v_1, v_2,v_N\}]
\end{align}
This open book simplicial complex has $N-2$ triangular pages $\{v_1, v_2,v_j\}$ with $j=3,\ldots,N$ and the bonding edge $\{v_1, v_2\}$. As an illustration, we show on Fig.~\ref{fig:book_illust} the open book with four triangular pages. 

For the RRSC(3)
\begin{equation}
\label{SN123:bounds}
\begin{split}
N-1 & \leq S_1 \leq 3N-6 \\
   0  & \leq S_2 \leq 3N-8 \\
  0    &\leq S_3 \leq N-3 \\
2N-1&\leq S \leq 8N-17
\end{split}
\end{equation}
The lower bound is achieved on trees. The upper bound is achieved on the simplicial complexes with maximal number of tetrahedrons, e.g., an open book complex
\begin{align}
\label{book-3}
 [\{v_1, v_2,v_3,v_4\}, \{v_1, v_2,v_3,v_5\}, \ldots, \{v_1, v_2,v_3,v_N\}]
\end{align}
This open book simplicial complex has $N-3$ tetrahedronal pages $\{v_1, v_2,v_3,v_j\}$ with $j=4,\ldots,N$ glued by the bonding triangle $\{v_1, v_2, v_3\}$. 

Generally for the RRSC(m), one can derive bounds similar to the bounds \eqref{SN12:bounds} and \eqref{SN123:bounds}. We mention only the bounds 
\begin{equation}
\label{SN:bounds-d}
2N-1\leq S \leq 2^m N- [(m-1)2^m +1]
\end{equation}
for the total number of simplices. The lower bound again occurs on trees. The upper bound occurs e.g. on the $m-$dimensional open book simplicial complex
\begin{equation}
\label{open-m}
[\{v_1,\ldots,v_m,v_j\},\, j=m+1,\ldots,N]
\end{equation}
with $m-$dimensional pages $\{v_1,\ldots,v_m,v_j\}$ glued by the $(m-1)-$dimensional simplex $\{v_1,\ldots,v_m\}$. 

The random quantities $S$ and $S_d$ are asymptotically self-averaging.
\rev{Formally, a random variable is said to be asymptotically self-averaging if its relative standard deviation, $\sqrt{\text{Var}[S_d]} / \left<S_d \right>$, approaches zero as $N \to \infty$. The most common scaling scaling}
\begin{equation} 
\label{eq:selfaveraging}
\frac{\sqrt{\text{Var}[S_d]}}{\left<S_d \right>} \sim N^{-1/2} 
\end{equation} 
\revv{holds for the RRSC($m$) models when $m\le5$, while for $m>5$, the relative standard deviation vanishes more slowly than $N^{-1/2}$; see Sec.~\ref{sec:urn}.} The asymptotic self-averaging implies that the average values $S(N)=\langle S\rangle$ and $S_d(N)=\langle S_d\rangle$ provide the chief information about the random quantities $S$ and $S_d$. The above bounds \eqref{SN12:bounds}, \eqref{SN123:bounds} and \eqref{SN:bounds-d} also suggest that the average values grow linearly with $N$:
\begin{subequations}
\begin{equation}
\label{SN-av}
S(N) \simeq AN
\end{equation}
and 
\begin{equation}
\label{SN-d-av}
S_d(N) \simeq a_d N
\end{equation}
Using $S_0=N$ and \eqref{SN:def} we get
\begin{equation}
\label{Aa-0}
a_0=1, \quad  A = \sum_{d=0}^m a_d
\end{equation}
The Euler characteristic remains constant, Eq.~\eqref{Euler:def}, and hence the alternating sum of the amplitudes vanishes
\begin{equation}
\label{Euler:amp}
\sum_{d=0}^m (-1)^d a_d = 0
\end{equation}
\end{subequations}

In Sec.~\ref{sec:linear}, we provide strong support for the linear growth laws \eqref{SN-av} and \eqref{SN-d-av} in all models RRSC(m) and show how to compute $a_d$ and $A$.

\section{Average numbers of simplices}
\label{sec:linear}

In this section, we consider the RRSC(m) models with $m\geq 2$. We shall assume asymptotic self-averaging of the total number of simplices of various dimensions, $S_1, \ldots, S_m$, and establish the linear growth of the average values $S_d(N)=\langle S_d\rangle$. \revv{This assumption is justified rigorously in Sec.~\ref{sec:urn}.} 

For the RRSC(1), the number of edges $S_1$ is strongly self-averaging, viz., deterministic:  $S_1=N-1$ since the RRSC(1) process generates trees. The exact results \eqref{S01} imply 
\begin{equation}
A=2,\quad a_0=a_1=1, \qquad a_k=0 \quad (k\geq 2)
\end{equation}

The RRSC(m) models with $m\geq 2$ are more challenging because the quantities $S_d$ with $d=1,\ldots, m-1$ are random for sufficiently large $N$. 

\subsection{RRSC(2)}

We have $S_0 = N$, while $S_1$ and $S_2$ are random. Adding a new vertex, $N\to N+1$, leads to 
\begin{equation}
\label{S12:rec}
\begin{split}
&(S_2,S_1)\to(S_2,S_1+1)       \quad\qquad\text{prob} \quad \frac{N}{N+S_1}\\
& (S_2,S_1)\to(S_2+1,S_1+2)  ~\quad\text{prob} \quad \frac{S_1}{N+S_1}
\end{split}
\end{equation}
Averaging \eqref{S12:rec} we obtain
\begin{equation}
\label{SN12:exact}
\begin{split}
&S_2(N+1)= S_2(N)+\left\langle\frac{S_1}{N+S_1}\right\rangle \\
&S_1(N+1)= S_1(N)+\left\langle\frac{N+2S_1}{N+S_1}\right\rangle
\end{split}
\end{equation}
Massaging Eqs.~\eqref{SN12:exact} yields
\begin{subequations}
\begin{align}
\label{SN2:exact}
& S_2(N+1) = S_2(N) + 1 - \left\langle\frac{N}{N+S_1}\right\rangle \\
\label{SN1:exact}
& S_1(N+1) = S_1(N) + 2 - \left\langle\frac{N}{N+S_1}\right\rangle  
\end{align}
\end{subequations}
Subtracting \eqref{SN2:exact} from \eqref{SN1:exact} we deduce the recurrence
\begin{equation}
S_1(N+1) - S_2(N+1) = S_1(N) -  S_2(N) + 1
\end{equation}
which has a simple solution
\begin{equation}
\label{SN12:Euler}
S_1(N) -  S_2(N) = N - 1
\end{equation}
To appreciate \eqref{SN12:Euler}, we note that reducing the Euler characteristic \eqref{Euler:def} to $m=2$ gives a relation $S_1-S_2=N-1$ between random quantities $S_1$ and $S_2$; averaging this relation leads to  \eqref{SN12:Euler}.

\begin{figure}
	\centering
	\includegraphics[width=1\linewidth]{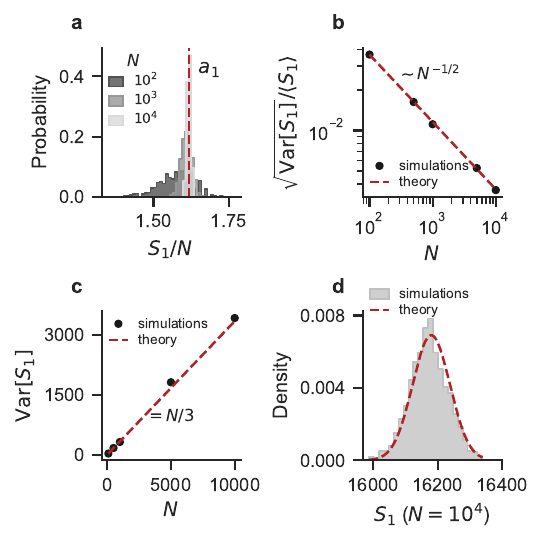}
\caption{Concentration and asymptotic self-averaging of $S_1$ for the RRSC(2) model. Black dots (b,~c) and gray histograms (a,~d) are simulations, each computed from 1000 independent RRSC(2) realizations; red dashed lines are theory.
		(a)~Distributions of $S_1/N$ at $N=10^2$, $10^3$, and $10^4$ concentrate sharply around $a_1=(1+\sqrt{5})/2$ (\cref{a1:2}) as $N$ grows. 
		(b)~The relative fluctuation $\sqrt{\text{Var}[S_1]}/\langle S_1\rangle$ decays as $\sim N^{-1/2}$, 
		the signature of the most common asymptotic self-averaging, see \eqref{eq:selfaveraging}.
		(c)~$\mathrm{Var}[S_1]=N/3$  (\cref{S1:var}).
		(d)~Distribution of $S_1$ at $N=10^4$ against the predicted Gaussian with mean $a_1 N$ and variance $N/3$ given by \cref{PSN:Gauss}.
	}
	\label{fig:S1}
\end{figure}

Equation \eqref{SN1:exact} involves only $S_1$, but it is not a closed recurrence for $S_1(N)=\langle S_1\rangle$ as we cannot express the average $\left\langle\frac{N}{N+S_1}\right\rangle$ through $S_1(N)$. Fortunately, the random variable $S_1$ is asymptotically self-averaging, as we show analytically in Appendix \ref{ap:edges}. \rev{We also confirm this numerically by simulating RRSC(2) structures for a range of $N$ values, with 1000 independent realizations for each $N$ (see \cref{fig:S1}a,b). Note that this is the procedure used throughout the text for numerical simulations.} Therefore $\left\langle\frac{N}{N+S_1}\right\rangle$ simplifies to $\frac{N}{N +S_1(N)}$ and Eq.~\eqref{SN1:exact} becomes
\begin{equation}
\label{SN1:eq}
S_1(N+1)= S_1(N)+ 2 - \frac{N}{N +S_1(N)} 
\end{equation}
when  $N\gg 1$. For large $N$, we can employ the continuum approach and reduce \eqref{SN1:eq} to an ordinary differential equation (ODE)
\begin{equation}
\label{SN1:ODE}
\frac{dS_1}{dN}=2-\frac{N}{N +S_1}
\end{equation}
which has indeed a linearly growing with $N$ solution, $S_1(N)=a_1 N$, with amplitude satisfying 
\begin{equation}
a_1= \frac{1+2a_1}{1+a_1}
\end{equation}
Thus 
\begin{equation}
\label{a1:2}
a_1 = \frac{\sqrt{5} +1}{2}\,, \quad a_2 = \frac{\sqrt{5}-1}{2}\,, \quad A=\sqrt{5} +1
\end{equation}
\rev{which we confirmed numerically (see \cref{fig:S1}a for $a_1$ \revv{and Appendix~\ref{ap:edges} for the other amplitudes})}.
Of course, $a_0=1$ and $a_k=0$ for $k\geq 3$.

In Appendix \ref{ap:edges}, we derived the variance of $S_1$: 
\begin{equation}
\label{S1:var}
\text{Var}[S_1] = \frac{1}{3}\, N
\end{equation}
\rev{which matches numerical simulations (see \cref{fig:S1}c,d)}.

\subsection{RRSC(3)}

We have $S_0 = N$, while $S_1, S_2, S_3$ are random. Adding a new vertex, $N\to N+1$, leads to
\begin{equation*}
(S_3,S_2,S_1)\to(S_3,S_2,S_1+1) 
\end{equation*}
with probability $\frac{S_0}{S-S_3}$,
\begin{equation*}
(S_3,S_2,S_1)\to(S_3,S_2+1,S_1+2)
\end{equation*}
with probability $\frac{S_1}{S-S_3}$, and
\begin{equation*}
(S_3,S_2,S_1)\to(S_3+1,S_2+3,S_1+3)
\end{equation*}
with probability $\frac{S_2}{S-S_3}$. Relying on the asymptotic self-averaging of $S_1$ and $S_2$ established in Appendix \ref{ap:edges} we arrive at 
\begin{equation}
\label{SN123:eq}
\begin{split}
&S_3(N+1)= S_3(N)+\frac{S_2(N)}{N+S_1(N)+S_2(N)} \\
&S_2(N+1)= S_2(N)+ \frac{S_1(N)+3S_2(N)}{N+S_1(N)+S_2(N)} \\
&S_1(N+1)= S_1(N)+ \frac{S_0(N)+2S_1(N)+3S_2(N)}{N+S_1(N)+S_2(N)}
\end{split}
\end{equation}
The last two equations lead to the algebraic equations for the amplitudes $a_1$ and $a_2$
\begin{equation}
\begin{split}
\label{a12}
&a_2= \frac{a_1+3a_2}{1+a_1+a_2}\\
&a_1= \frac{1+2a_1+3a_2}{1+a_1+a_2}
\end{split}
\end{equation}
and the first equation in \eqref{SN123:eq} expresses $a_3$ via $a_1$ and $a_2$:
\begin{align}
\label{a3}
a_3= \frac{a_2}{1+a_1+a_2}
\end{align}
Equations \eqref{a12} constitute a closed system which does not seem solvable in quadratures. One can compute $a_1$ and $a_2$ to any desired accuracy, however, and then use \eqref{a3} to determine $a_3$. Thus one gets 
\begin{equation}
\label{aA:3}
\begin{split}
&a_1 = 2.089\,396\,431\ldots\\
&a_2 = 1.401\,467\,183\ldots\\
&a_3 = 0.312\,070\,751\ldots\\
&A = 4.802\,934\,366\ldots
\end{split}
\end{equation}
where $A=1+a_1+a_2+a_3$, cf. \eqref{Aa-0}. 
\revv{Numerical simulations confirm these values (Appendix~\ref{ap:edges}).}

\subsection{RRSC(4)}

The same analysis as before yields a closed system of algebraic equations for the first three amplitudes:
\begin{equation}
\label{a123}
\begin{split}
&a_3= \frac{a_2+4a_3}{1+a_1+a_2+a_3} \\
&a_2= \frac{a_1 + 3a_2+6a_3}{1+a_1+a_2+a_3} \\
&a_1= \frac{1+2a_1+3a_2+4a_3}{1+a_1+a_2+a_3}
\end{split}
\end{equation}
The last amplitude $a_4$ is then given by 
\begin{equation}
\label{a4}
a_4= \frac{a_3}{1+a_1+a_2+a_3}
\end{equation}
One can find the amplitudes to any desired accuracy:
\begin{equation}
\label{aA:4}
\begin{split}
&a_1 = 2.450\,304\,492\ldots\\
&a_2 = 2.188\,116\,583\ldots\\
&a_3 = 0.871\,713\,023\ldots\\
&a_4 = 0.133\,900\,931\ldots\\
&A = 6.644\,035\,030\ldots
\end{split}
\end{equation}
\revv{Numerical simulations confirm these values (Appendix~\ref{ap:edges}).} \rev{Note that for larger $m$ convergence is slower, i.e. require larger $N$, because building larger simplices required the presence of simplices of all of the smaller dimensions.}

Extending the analysis to the RRSC(m) models with $m\geq 5$ is straightforward. The amplitudes in the models with $m\leq 10$ are listed in 
Appendix \ref{ap:10}.

\subsection{General case}

For the RRSC(m) model, one obtains a closed system of $m-1$ algebraic equations
\begin{equation}
\label{ak-m}
a_k \left[1+\sum_{j=1}^{m-1}a_j\right] =\sum_{d=k-1}^{m-1}a_d\binom{d+1}{k}
\end{equation}
for $k=1,\ldots, m-1$. This system of equations reduces to Eqs.~\eqref{a12} when $m=3$ and Eqs.~\eqref{a123} when $m=4$. The last amplitude is found from [cf. \eqref{a3} and \eqref{a4}]
\begin{equation}
\label{am-m}
a_m = a_{m-1}\left[1+\sum_{j=1}^{m-1}a_j\right]^{-1}
\end{equation}
As always, we have $a_0=1$ and $A = \sum_{0\leq j\leq m}a_j$.

\revv{Although Eqs.~\eqref{ak-m} look nonlinear, they are equivalent to a linear eigenvalue problem, as we show next.}

\subsection{\texorpdfstring{\revv{Connection to generalized P\'olya urns}}{Connection to generalized Polya urns}}
\label{sec:urn}

\begin{figure}[t]
	\centering
	\includegraphics[width=1\linewidth]{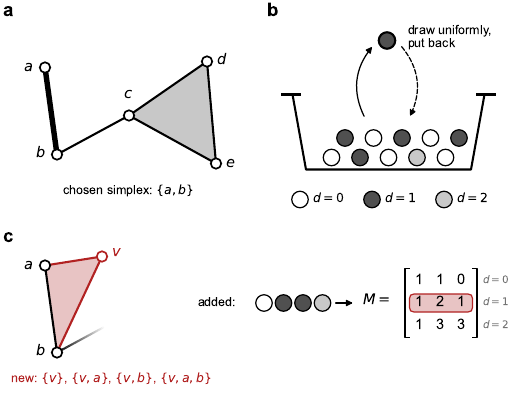}
	\caption{\revv{Mapping of the RRSC($m$) model onto a generalized P\'olya urn, illustrated for $m=3$. (a) A simplicial complex; a new vertex $v$ joins the chosen simplex $\{a,b\}$. (b) Each simplex of dimension $d<m$ is a ball of type $d$. Choosing a simplex uniformly amounts to drawing a ball uniformly and putting it back. (c) The new vertex $v$ creates the simplices $\{v\}\cup\tau$ for all $\tau\subseteq\{a,b\}$: one vertex, two edges, and one triangle. These new balls form row $d=1$ of the replacement matrix~(\ref{eq:urnM}).}}
	\label{fig:urnscheme}
\end{figure}

\revv{The self-averaging assumption used above can be removed, because the RRSC($m$) process is a generalized P\'olya urn~\cite{AthreyaKarlin68,Janson04,Mahmoud08}. Regard each simplex of dimension $d<m$ as a ball of type $d$ (Fig.~\ref{fig:urnscheme}). Choosing a simplex uniformly then amounts to drawing a ball uniformly. Since every new simplex contains the new vertex, drawing a $d$-simplex adds exactly $\binom{d+1}{k}$ new $k$-simplices, so the urn has the deterministic replacement matrix}
\begin{equation}
\revv{M_{dk}=\binom{d+1}{k},\qquad d,k=0,\ldots,m-1.}
\label{eq:urnM}
\end{equation}
\revv{This matrix is non-negative and irreducible, and limit theorems for such urns~\cite{Janson04} yield three results valid for every $m$. First, $S_d/N\to a_d$ almost surely, so that all $S_d$ are extensive and self-averaging. Second, the amplitudes form the Perron eigenvector of $M$:}
\begin{equation}
\revv{\sum_{d=0}^{m-1}a_dM_{dk}=\lambda_1a_k,\quad \lambda_1=1+a_1+\cdots+a_{m-1}}
\label{eq:perron}
\end{equation}
\revv{where $\lambda_1$ is the largest eigenvalue of $M$. Equations~(\ref{ak-m}) thus become a linear eigenvalue problem with a unique positive solution that reproduces the amplitudes of Appendix~\ref{ap:10}. Third, the fluctuations are controlled by the ratio $\lambda_2/\lambda_1$ of the two largest eigenvalues, listed in Table~\ref{tab:ratio}. When $\lambda_2/\lambda_1<1/2$, i.e., for $m\le5$, the fluctuations of $S_d$ are Gaussian on the scale $\sqrt{N}$. For $m\ge6$, $S_d-a_dN$ grows as $N^{\lambda_2/\lambda_1}$ and has a non-Gaussian limit law; relative fluctuations still vanish, but more slowly than $N^{-1/2}$. In Appendix~\ref{app:urn}, we present the derivation, an intuitive explanation of the threshold $1/2$, and numerical checks.}

\begin{table}
\caption{\revv{Largest eigenvalue $\lambda_1$ of the matrix~(\ref{eq:urnM}) and ratio $\lambda_2/\lambda_1$ of the two largest eigenvalues, for the RRSC($m$) models.}}
\label{tab:ratio}
\revv{%
\begin{tabular*}{\columnwidth}{@{\extracolsep{\fill}}ccccccccc}
\toprule
$m$ & 2 & 3 & 4 & 5 & 6 & 7 & 8 & 9\\
\midrule
$\lambda_1$ & 2.618 & 4.491 & 6.510 & 8.626 & 10.812 & 13.050 & 15.329 & 17.643\\
$\lambda_2/\lambda_1$ & 0.146 & 0.299 & 0.407 & 0.485 & 0.542 & 0.587 & 0.622 & 0.651\\
\bottomrule
\end{tabular*}}
\end{table}

\section{Primordial vertex}
\label{sec:primordial}

The evolution begins with a single primordial vertex $v_1$. This vertex plays a prominent role because it is present throughout the evolution. The second vertex $v_2$ plays a comparable role since, in the simplicial complex with two vertices, $\mathcal{K}_2=\langle\{v_1\}, \{v_2\}; \{v_1, v_2\}\rangle$, we can relabel the vertices: $v_1\leftrightarrow v_2$. For concreteness, we consider characteristics involving the primordial vertex. 

Denote by $\delta_1$ the degree of the primordial vertex, i.e., the number of edges of type $\{v_1,v_k\}$ in $S_1$. The degree $\delta_1$ is a random quantity varying in the range
\begin{equation}
\label{delta:extreme}
1\leq \delta_1\leq N-1
\end{equation}
when $N\geq 2$. Indeed, the $\{v_1,v_2\}$ edge is present for all $N\geq 2$. The upper bound is achieved, e.g., on the star tree \eqref{star}. One would like to compute the probabilities to observe the extremal values in \eqref{delta:extreme}. The most important and interesting challenge is to determine the probability distribution $P(\delta_1,N)$.

Denote by $\delta_j$ the degree of vertex $v_j$, equivalently the number of edges $\{v_j,v_k\}$  in $S_1$. The degree $\delta_1$ is, on average, higher than $\delta_j$ with $2<j\leq N$. The degrees $\delta_1$ and $\delta_2$ have the same distributions for $N\geq 2$. One can ask about the probability that the degree of the primordial vertex is the largest
\begin{equation}
\label{FN:def}
F(N)=\text{Prob}[\delta_1\geq \delta_j, ~ j=2,\ldots,N]
\end{equation}
at the end of the evolution. Another natural quantity is the probability that the degree of the primordial vertex was the largest throughout the evolution:
\begin{equation}
\Phi(N)=\text{Prob}[\delta_1\geq \delta_j, ~ j=2,\ldots, n \,|\, n=3,\ldots,N]
\end{equation}
Even for the RRSC(1), i.e., the RRT, we only know \cite{KR02} the probability distribution $P(\delta_1,N)$ and the asymptotic behavior of $F(N)$; the `survival' probability $\Phi(N)$ is still unknown. 

Although we do not know the probability distribution $P(\delta_1,N)$ for the  RRSC(m) models with $m\geq 2$, it is possible to establish the growth law of the average degree $\delta_1(N)=\langle \delta_1\rangle$ of the primordial vertex. We then assume that the average degree of the primordial vertex is comparable to the maximal degree. Suppose the degree distribution has an algebraic tail:
\begin{equation}
\label{deg:tail}
n_\delta\sim \delta^{-\nu}\qquad\text{when}\quad \delta\gg 1
\end{equation}
The maximal degree is estimated from the criterion
\begin{equation}
\label{deg:max-crit}
N\sum_{\delta\geq \delta_\text{max}}n_\delta \sim 1
\end{equation}
which in conjunction with \eqref{deg:tail} give
\begin{equation}
\label{deg:max}
\delta_\text{max} \sim N^\frac{1}{\nu -1}
\end{equation}
Assuming that the average degree of the primordial vertex has the same scaling with $N$ as the maximal allows one to fix the decay exponent $\nu$. 

The computations for the RRSC(m) models quickly become cumbersome when $m$ increases, as one needs to know the amplitudes $a_d(m)$ with $d=1,\ldots,m-1$ [cf. Sec.~\ref{sec:linear}] and the maximal eigenvalue of an $(m-1)\times (m-1)$ matrix. \revv{In Sec.~\ref{sec:nu-general}, we show that the result for general $m$ reduces to the ratio of Perron eigenvalues of two integer matrices.} We also emphasize that our approach only gives the tail of the distribution.  

If the degree distribution has an algebraic tail \eqref{deg:tail}, the exponent must obey $\nu>2$. To establish this inequality, we mention that the sum of degrees is twice the number of edges:
\begin{equation}
\label{sum}
\sum_{j=1}^N \delta_j = 2S_1(N)
\end{equation}
 Rewriting the normalization condition and \eqref{sum} via the degree distribution $n_\delta$ give the sum rules 
\begin{equation}
\label{sum:rules}
\sum_{\delta\geq 1}n_\delta = 1, \qquad  \sum_{\delta\geq 1}\delta n_\delta = 2 a_1
\end{equation}
The inequality $\nu>2$ is needed to ensure that the second sum in \eqref{sum:rules} is finite.

\subsection{RRSC(1)}
\label{subsec:1} 

Here, we present some results about the evolution of the degrees $\delta_1$ and $\delta_\text{max}$ in the RRSC(1) model. These results are mostly known \cite{KR02}. We remind them to gain intuition about the nature of the random variables $\delta_1$ and $\delta_\text{max}$ in the tractable model and guess the conjectural behaviors of these random variables in the RRSC(m) models with $m\geq 2$. 

Adding a new vertex, $N\to N+1$, leads to the stochastic recurrence equation 
\begin{equation}
\label{delta-RRT}
\delta_1\to 
\begin{cases}
\delta_1+1 & \text{prob} \quad \frac{1}{N}\\
\delta_1 & \text{prob} \quad 1-\frac{1}{N} 
\end{cases}
\end{equation}

Using Eq.~\eqref{delta-RRT}, we express the cumulant generating function
\begin{equation}
\label{CG}
\left\langle e^{\lambda \delta_1}\right\rangle = \sum_{\delta_1=1}^{N-1} e^{\lambda \delta_1} P(\delta_1,N)
\end{equation}
via the ratio of gamma functions
\begin{equation}
\label{CG:sol}
\left\langle e^{\lambda \delta_1}\right\rangle = \frac{\Gamma(e^\lambda - 1 +N)}{\Gamma(e^\lambda + 1)\,\Gamma(N)}
\end{equation}
from which one extracts
\begin{equation}
\label{P:Stirling}
P(\delta_1,N)=\frac{1}{(N-1)!}\,{N-1\brack \delta_1}
\end{equation}
Here  ${n\brack m}$ are the Stirling numbers of the first kind \cite{Knuth}. Recalling the identities ${n\brack 1}=(n-1)!$ and  ${n\brack n}=1$, we find that the extreme values \eqref{delta:extreme} are achieved with probabilities
\begin{subequations}
\label{min-max:RRT}
\begin{align}
\label{min:primordial}
&P(1,N)=\frac{1}{N-1}\\
\label{max:primordial}
& P(N-1,N)=\frac{1}{(N-1)!}
\end{align}
\end{subequations}

The standard relation
\begin{equation}
\label{CG:exp}
\ln \left\langle e^{\lambda \delta_1}\right\rangle  = \sum_{n\geq 1} \frac{\lambda^n}{n!}\,  \langle\!\langle \delta_1^n\rangle\!\rangle
\end{equation}
gives the cumulants: The average $\delta_1(N)=\langle \delta_1\rangle$, the variance $V_1(N)=\text{Var}[\delta_1]=\langle\!\langle \delta_1^2\rangle\!\rangle\equiv \langle \delta_1^2\rangle-\langle \delta_1\rangle^2$, and higher cumulants $\langle\!\langle \delta_1^n\rangle\!\rangle$ with $n\geq 3$. The first two cumulants admit exact expressions (valid for all $N$) 
\begin{equation}
\delta_1(N) = H_{N-1}, \quad  V_1(N)=H_{N-1} - H_{N-1}^{(2)} 
\end{equation}
where $H_k=\sum_{1\leq j\leq k}j^{-1}$ and $H_k^{(2)}=\sum_{1\leq j\leq k}j^{-2}$ are harmonic numbers and $\gamma$ is the Euler constant. The asymptotic behaviors are
\begin{equation}
\begin{split}
&\delta_1(N) = \ln N +\gamma +O(N^{-1})\\
&V_1(N) = \ln N +\gamma  -\frac{\pi^2}{6}+O(N^{-1})
\end{split}
\end{equation}

Using \eqref{CG:sol} and \eqref{CG:exp} one can extract cumulants $\langle\!\langle \delta_1^n\rangle\!\rangle$ for all $n\geq 1$. The exact answers contain generalized harmonic numbers. The leading behaviors are remarkably universal: 
\begin{equation}
\label{cum}
\langle\!\langle \delta_1^n\rangle\!\rangle \simeq \ln N 
\end{equation}
We thus conclude that the random variable $\delta_1$ is asymptotically self-averaging. Furthermore, the growth laws \eqref{cum} imply that the distribution $P(\delta_1,N)$ is asymptotically Gaussian
\begin{equation}
\label{P:Stirling-Gauss}
P(\delta_1,N) \simeq \frac{1}{\sqrt{2\pi\ln N}}\,\exp\!\left[-\frac{(\delta_1-\ln N)^2}{2\ln N}\right]
\end{equation}
when $N\gg 1$. Alternatively, the Gaussian form \eqref{P:Stirling-Gauss} can be directly deduced from Eq.~\eqref{delta-RRT}.

The degree distribution for the RRSC(1) is $2^{-\delta}$ which, in conjunction with criterion \eqref{deg:max}, gives $\delta_\text{max} \simeq \log_2 N$. Thus, for the RRT, the maximal degree is comparable with the degree of the primordial vertex. The distribution $\Pi(\delta_\text{max}, N)$ is very different from the asymptotically Gaussian distribution \eqref{P:Stirling-Gauss} of the degree of the primordial vertex. The variance $\text{Var}[\delta_\text{max}]=\langle\!\langle  \delta_\text{max}^2\rangle\!\rangle$ is unknown, but expected to remain finite in the $N\to\infty$ limit.  

When $N\geq 3$, the maximal degree varies in the range $2\leq \delta_\text{max}\leq N-1$. The extreme values \eqref{delta:extreme} are achieved with probabilities
\begin{equation}
\label{extreme:max}
\Pi(2,N)=\frac{2^{N-2}}{(N-1)!}\,, \quad \Pi(N-1,N)=\frac{2}{(N-1)!}
\end{equation}
The maximal possible degree $\delta_\text{max}=N-1$ is realized when either $\deg(v_1)=N-1$ or $\deg(v_2)=N-1$. These two outcomes are equally plausible, and only one may occur when $N\geq 3$. Therefore, the probability twice exceeds the corresponding probability \eqref{max:primordial}. The minimal value $\delta_\text{max}=2$ of the maximal degree occurs when the emerging tree is a linear graph. The probability of this event is found from the recurrence 
\begin{equation}
\Pi(2,N+1)=\frac{2}{N}\,\Pi(2,N)
\end{equation}
and initial condition $P(2,3)=1$. 

The probability $\Pi(2,N)$ of the minimal value of the maximal degree is much smaller than the corresponding probability \eqref{min:primordial}. This feature is a quantitative sign of the sharpness of the distribution $\Pi(\delta_\text{max}, N)$ of the maximal degree. In comparison, the distribution $P(\delta_1, N)$ is much less sharp. 

The sharpness of the distribution $\Pi(\delta_\text{max},N)$ allows the computation of the asymptotic behavior of the probability $F(N)$ that the degree of the primordial vertex is the largest at the end of the evolution:
\begin{equation}
F(N)\simeq \sum_{\delta_1> \log_2 N}P(\delta_1,N)
\end{equation}
Using \eqref{P:Stirling} and the properties of the Stirling numbers of the first kind \cite{Knuth} one finds  \cite{KR02}
\begin{equation}
\label{FN:RRT}
F(N) \sim N^{-\epsilon}(\ln N)^{-\frac{1}{2}}
\end{equation}
with exponent 
\begin{equation}
\label{eps}
\epsilon=1-\frac{1 + \ln(\ln 2)}{\ln 2}= 0.086071332\ldots
\end{equation}
appearing in a surprisingly large number of unrelated problems, e.g., in number theory \cite{Ford08}, averaging processes \cite{BK21}, and random permutations \cite{Peres16,Ford16,Ford22}.

\subsection{RRSC(2)}
\label{subsec:2}

Adding a new vertex, $N\to N+1$, leads to the stochastic recurrence equation 
\begin{equation}
\label{delta-change}
\delta_1\to 
\begin{cases}
\delta_1+1 & \text{prob} \quad \frac{1+\delta_1}{N+S_1}\\
\delta_1 & \text{prob} \quad 1-\frac{1+\delta_1}{N+S_1}
\end{cases}
\end{equation}
for the degree of the primordial vertex. Indeed, the arriving vertex $v_{N+1}$ may join the primordial vertex to form the edge $\{v_1,v_{N+1}\}$. The probability of this event is $1/[N+S_1]$. The probability $\delta_1/[N+S_1]$ accounts for $v_{N+1}$ joining an edge of the type $\{v_1,v_k\}$ in $S_1$. Thanks to these two events, the degree of the primordial vertex increases by one. 

The stochastic equation \eqref{delta-change} is much more complicated than the corresponding stochastic equation \eqref{delta-RRT} for the RRSC(1) model. Indeed, the probabilities in Eq.~\eqref{delta-RRT} depend only on $N$; the probabilities in Eq.~\eqref{delta-change} depend on $N+S_1$, with $S_1$ determined by the evolution history. Thus, the evolution process \eqref{delta-change} is non-Markovian, making analytical advances much more challenging. For instance, iterating Eq.~\eqref{delta-RRT} one readily derives \eqref{min-max:RRT}. In contrast, for the RRSC(2) model, the extreme probabilities $P(1, N)$ and $P(N-1, N)$ are unknown. 
 
Using Eq.~\eqref{delta-change} we find that the average degree of the primordial vertex increases according to
\begin{equation}
\label{delta-rec}
\delta_1(N+1) = \delta_1(N)+ \left\langle \frac{1+\delta_1}{N+S_1}\right\rangle
\end{equation}
Using the asymptotic self-averaging of $S_1$ we re-write Eq.~\eqref{delta-rec} as
\begin{equation}
\label{delta-recur}
\delta_1(N+1) = \delta_1(N)+ \frac{1+\delta_1(N)}{N+S_1(N)}
\end{equation}
The replacement of Eq.~\eqref{delta-rec} by Eq.~\eqref{delta-recur} ignores correlations between the random quantities $\delta_1$ and $S_1$. As an a posteriori justification, we note that using Eq.~\eqref{delta-recur} we establish a sub-linear growth of $\delta_1$, see \eqref{deg1}. Thus, the influence of correlations of $\delta_1$ and a linearly growing with $N$ asymptotically self-averaging random quantity $S_1$ is asymptotically negligible. 

\begin{figure}[h]
	\centering
	\includegraphics[width=1\linewidth]{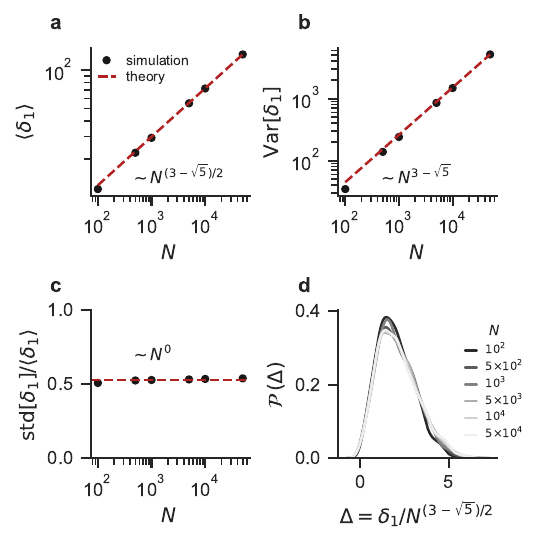}
	\caption{Distribution of the degree $\delta_1$ of the primordial vertex, as $N$ grows, for RRSC(2). 
	Black dots (a, b,~c) and gray distributions (d) are simulations, each computed from 1000 independent RRSC(2) realizations; red dashed lines are theory. 
		(a)~The average primordial degree scales as $\langle\delta_1\rangle\sim
		N^{(3-\sqrt{5})/2}$ (\cref{deg1}).
		(b) Its variance scales as $\mathrm{Var}[\delta_1] \sim N^{3-\sqrt{5}}$ (\cref{deg1-var}); the exponent is exactly twice that of the average.
		(c)~As a consequence, the relative fluctuation is constant in $N$, indicating that $\langle\delta_1\rangle$ is not self-averaging. 
		(d)~Collapse of the rescaled distributions of the primordial degree onto $\mathcal{P}(\Delta)$ (\cref{eq:rescaled_distribution}).
	}
	\label{fig:primordial_degree}
\end{figure}

Equation \eqref{delta-recur} is applicable when $N\to \infty$. In this limit, we can adopt a continuum description and transform Eq.~\eqref{delta-recur} into an ODE
\begin{equation}
\label{delta-eq}
\frac{d \delta_1}{d N}=\frac{1+\delta_1}{(1+a_1)N}
\end{equation}
Integrating \eqref{delta-eq} and recalling that $a_1=\frac{1+\sqrt{5}}{2}$ for the RRSC(2), see \eqref{a1:2}, we find that the average degree of the primordial vertex increases algebraically with size:
\begin{equation}
\label{deg1}
\delta_1(N) \sim N^\frac{3-\sqrt{5}}{2}
\end{equation}
\rev{which we confirmed numerically (\cref{fig:primordial_degree}a)}.

\begin{figure*}[htb]
	\centering
	\includegraphics[width=1\linewidth]{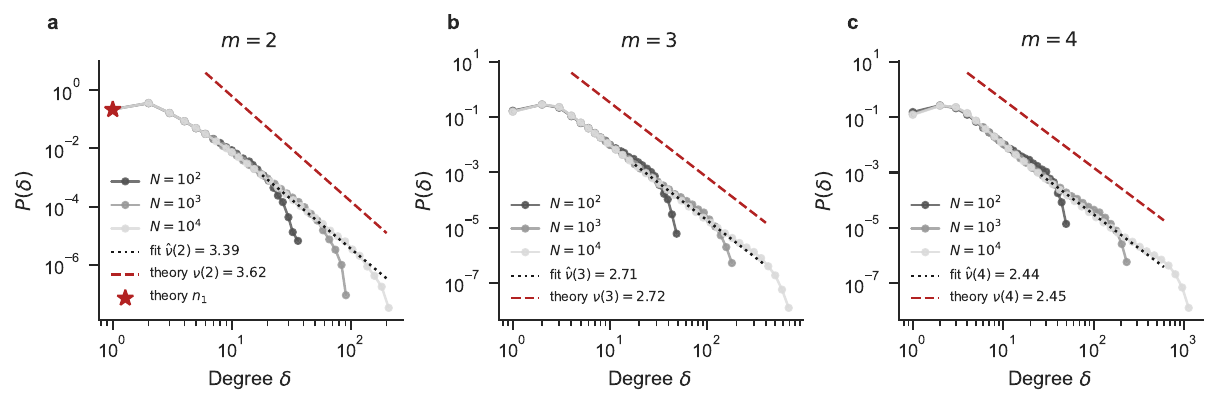}
	\caption{Degree distribution for (a) RRSC(2), (b) RRSC(3), and (c) RRSC(4). 
	In each case, degree distributions are computed from numerical simulations for $N=10^2$, $10^3$, and $10^4$ (shades of gray), with 1000 independent realizations for each $N$. Numerical exponents $\hat \nu(m)$ are computed from fitting the data (dotted black line) and theoretical scalings (dashed red line) are obtained from the theoretical exponents in \cref{tail-2,tail-3,tail-4}, respectively.
	In panel (a), numerics also match the theoretical value of $n_1$ (\cref{deg:1-2}).}
	\label{fig:degree}
\end{figure*}

Assuming that the degree distribution has an algebraic tail \eqref{deg:tail}, we find that the maximal degree grows as $\delta_\text{max} \sim N^\frac{1}{\nu -1}$, see \eqref{deg:max-crit}--\eqref{deg:max}. By definition, $\delta_\text{max} \geq \delta_1$. The maximal degree is expected to be comparable with the degree of the primordial vertex $\delta_1$. Comparing \eqref{deg1} and \eqref{deg:max} we get $\frac{1}{\nu -1}=\frac{3-\sqrt{5}}{2}$ leading to the announced result \eqref{tail-2} for the exponent $\nu=\nu(2)$ in \eqref{deg:tail}.
\rev{We also computed this exponent by fitting the degree distributions obtained from numerical simulations for $N=10^2$, $10^3$, and $10^4$ (\cref{fig:degree}a). The numerical exponent $\hat \nu(2) \simeq 3.39$ which is close to the theoretical value $\nu(2)=\frac{5+\sqrt{5}}{2} \simeq 3.62$ although somewhat lower. This reflects the shorter power-law tail at $m=2$: the finite-$N$ cutoff biases the fit low. Indeed, at low $m$ each new vertex joins fewer existing vertices and creates fewer faces, so degrees grow slowly; and since a vertex belongs to fewer simplices, even the most-connected vertices gain edges rarely. This keeps the maximum degree, and hence the tail, short. 
As we will see below, the numerical exponents for larger $m=3$ and $m=4$ match almost exactly the theoretical values.}  

For the RRSC(1), the random variables $\delta_1$ and $\delta_\text{max}$ are asymptotically self-averaging, and the distribution of the maximal degree $\delta_\text{max}$ is very sharp (Sec.~\ref{subsec:1}). For the RRSC(2), and generally for the RRSC(m) models with $m\geq 2$,  the random variables $\delta_1$ and $\delta_\text{max}$ are non-self-averaging. The proof of this assertion is lacking, although the numerical evidence is quite compelling. 

To probe the behavior of the variance $V_1(N)=\text{Var}[\delta_1]$ we rely on the asymptotic self-averaging of $S$ and $S_2$ and employ the same arguments as in the derivation of the recurrence relation \eqref{delta-recur} for the average to yield  
\begin{equation}
\label{var-recur}
V_1(N+1) = \left[1+\frac{2}{(1+a_1)N}\right]V_1(N)+ \frac{3\delta_1(N)}{(1+a_1)N}
\end{equation}
where we have kept only two leading terms. An elementary analysis of \eqref{var-recur} shows that the second term on the right-hand side is also subdominant. Omitting the second term and adopting a continuum description, we simplify the recurrence \eqref{var-recur} to an ODE
\begin{equation}
\label{V1-eq}
\frac{d V_1}{d N}=\frac{2 V_1}{(1+a_1)N}
\end{equation}
from which
\begin{equation}
\label{deg1-var}
V_1(N) \sim N^{3-\sqrt{5}}
\end{equation}
\rev{which is also confirmed by numerical simulations (\cref{fig:primordial_degree}b)}.
Thus, the standard deviation $\sqrt{V_1(N)}$ scales similarly to the average indicating that $\delta_1$ is a non-self-averaging random variable \rev{(see the constant relative fluctuation in \cref{fig:primordial_degree}c)}. The probability distribution $P(\delta_1,N)$ is expected to approach the scaling form
\begin{equation}\label{eq:rescaled_distribution}
P(\delta_1,N) = N^{-\frac{3-\sqrt{5}}{2}}\,\mathcal{P}(\Delta)
\end{equation}
in the scaling limit
\begin{equation}
N\to\infty, \quad \delta_1\to\infty, \quad \Delta=\frac{\delta_1}{N^\frac{3-\sqrt{5}}{2}} = \text{finite}
\end{equation}
The scaled distribution $\mathcal{P}(\Delta)$ is \rev{not known analytically, but we obtain it numerically in \cref{fig:primordial_degree}d where distributions at different $N$ collapse onto a single scaling form.}

The probabilities $F(N)$ and $\Phi(N)$ presumably saturate, $F(\infty)>0$ and $\Phi(\infty)>0$, and hence not very interesting. These conjectural behaviors are suggested by similar behaviors observed in networks growing via preferential attachment \cite{KR02}. 

\subsection{RRSC(3)}
\label{subsec:3}

For the RRSC(3) model, to determine  the degree of the primordial vertex, we need to consider $\delta_1$ edges $\{v_1,v_j\}$ in $S_1$ as in the RRSC(2) model, and also $\Delta_1$ triangles $\{v_1,v_i,v_j\}$ in $S_2$.  The average values $\delta_1(N)= \langle \delta_1\rangle$ and $\Delta_1(N)= \langle \Delta_1\rangle$ satisfy the discrete evolution equations 
\begin{equation}
\label{delta-Delta-rec}
\begin{split}
&\delta_1(N+1) = \delta_1(N) + \left\langle \frac{1+\delta_1+\Delta_1}{S-S_3} \right\rangle\\
&\Delta_1(N+1) = \Delta_1(N) +  \left\langle \frac{\delta_1+2\Delta_1}{S-S_3} \right\rangle
\end{split}
\end{equation}
derived similarly to \eqref{delta-recur}. We proceed as before and assume the asymptotic self-averaging of $S_1,S_2,S_3$ and sub-linear growth of $\delta_1$ and $\Delta_1$ which we justify a posteriori. These assumptions allow us to reduce Eqs.~\eqref{delta-Delta-rec} to a pair of ODEs 
\begin{subequations}
\label{delta-Delta}
\begin{align}
\label{delta-3}
&\frac{d \delta_1}{d N}=\frac{\delta_1+\Delta_1}{(1+a_1+a_2)N}\\
\label{Delta-3}
&\frac{d \Delta_1}{d N}=\frac{\delta_1+2\Delta_1}{(1+a_1+a_2)N}
\end{align}
\end{subequations}
in the $N\to\infty$ limit. Solving \eqref{delta-Delta} we find that $\delta_1$ and $\Delta_1$ scale similarly 
\begin{equation}
\label{deg1-3}
\delta_1 \sim N^\frac{3+\sqrt{5}}{2(1+a_1+a_2)}, \qquad \Delta_1 \sim N^\frac{3+\sqrt{5}}{2(1+a_1+a_2)}
\end{equation}
The continuum approach does not allow one to fix the amplitudes in \eqref{deg1-3} as it was also the case in \eqref{deg1}. 
However, the ratio of $\Delta_1$ to $\delta_1$ saturates to the known number
\begin{equation}
\label{d-D:ratio}
\lim_{N\to\infty} \frac{\Delta_1(N)}{\delta_1(N)} = \frac{1+\sqrt{5}}{2}
\end{equation}
Indeed, dividing \eqref{Delta-3} by \eqref{delta-3} and treating $\Delta_1$ as a function of $\delta_1$ we get
\begin{equation}
\label{d-D:eq}
\frac{d\Delta_1}{d\delta_1} = \frac{\delta_1+2\Delta_1}{\delta_1+\Delta_1}
\end{equation}
from which one immediately deduces \eqref{d-D:ratio}. 

Repeating the same analysis as for the RRSC(2) we fix the exponent $\nu=\nu(3)$ in \eqref{deg:tail} for the RRSC(3):
\begin{equation}
\label{tail-3}
\nu(3)=1+ \frac{2(1+a_1+a_2)}{3+\sqrt{5}}=2.715\,357\,262\ldots
\end{equation}
\rev{matched by numerics in \cref{fig:degree}(b).}
[The amplitudes $a_1$ and $a_2$ are given by \eqref{aA:3}.]

\subsection{RRSC(4)}
\label{subsec:4}

For the RRSC(4) model, we denote by $\delta_1(N)$ the number of edges $\{v_1,v_j\}$ in $S_1(N)$, by $\Delta_1(N)$ the number of triangles $\{v_1,v_i,v_j\}$ in $S_2(N)$, and by $\tau_1(N)$ the number of tetrahedrons $\{v_1,v_i,v_j,v_k\}$ in $S_3(N)$. Similarly to \eqref{delta-Delta} we obtain
\begin{subequations}
\label{delta-Delta-tau}
\begin{align}
\label{delta-4}
&\frac{d \delta_1}{d N}=\frac{\delta_1+\Delta_1+\tau_1}{(1+a_1+a_2+a_3)N}\\
\label{Delta-4}
&\frac{d \Delta_1}{d N}=\frac{\delta_1+2\Delta_1+3\tau_1}{(1+a_1+a_2+a_3)N}\\
\label{tau-4}
&\frac{d \tau_1}{d N}=\frac{\Delta_1+3\tau_1}{(1+a_1+a_2+a_3)N}
\end{align}
\end{subequations}
in the $N\to\infty$ limit. From these equations we find the asymptotic behaviors
\begin{equation}
\label{deg1-4}
\delta_1\sim \Delta_1\sim \tau_1 \sim N^\frac{\Lambda}{1+a_1+a_2+a_3}
\end{equation}
and the saturated ratios
\begin{equation}
\label{two-ratios}
\begin{split}
&\lim_{N\to\infty} \frac{\Delta_1(N)}{\delta_1(N)} = \frac{\Lambda-3}{\Lambda^2-5\Lambda+3} = 2.089\,396\ldots  \\
&\lim_{N\to\infty} \frac{\tau_1(N)}{\delta_1(N)} = \frac{1}{\Lambda^2-5\Lambda+3} = 1.401\,467\ldots
\end{split}
\end{equation}
Here $\Lambda=4.490\,863\,615 \ldots$ is the largest eigenvalue of the matrix
\begin{equation}
\label{mat:3}
\begin{bmatrix}
    1 & 1 & 1 \\
    1 & 2 & 3 \\
    0 & 1 & 3
\end{bmatrix}
\end{equation}

The ratios \eqref{two-ratios} are derived from \eqref{delta-Delta-tau} similarly to the derivation of the ratio \eqref{d-D:ratio} in Sec.~\ref{subsec:3}. Namely, we divide \eqref{Delta-4} by \eqref{delta-4} and \eqref{tau-4} by \eqref{delta-4} and treat $\Delta_1$ and $\tau_1$ as functions of $\delta_1$. We find that $\Delta_1$ and $\tau_1$ scale linearly with $\delta_1$, and the ratios of  $\Delta_1$ and $\tau_1$ to $\delta_1$ are given by formulae \eqref{two-ratios} depending only on the largest eigenvalue of the matrix \eqref{mat:3}. 

Assuming again $\delta_1\sim \delta_\text{max}$ and comparing \eqref{deg:max} with \eqref{deg1-4} we fix the exponent $\nu=\nu(4)$ in \eqref{deg:tail}:
\begin{equation}
\label{tail-4}
\nu(4)=1+ \frac{1+a_1+a_2+a_3}{\Lambda}=2.449\,639\,680\ldots
\end{equation}
\rev{matched by numerics in \cref{fig:degree}(c).}
[The amplitudes $a_1,~a_2,~a_3$ are given by \eqref{aA:4}.]

\subsection{\texorpdfstring{\revv{General $m$}}{General m}}
\label{sec:nu-general}

\revv{The above results for $m\le4$ follow a general pattern valid for arbitrary $m$, but the derivations are rather pedestrian. A mapping to generalized P\'olya urns outlined in Sec.~\ref{sec:urn} allows one to consider arbitrary $m$. In Appendix~\ref{app:urn}, we provide more details. Namely, we show that the link of any vertex, i.e., the simplices that form a simplex together with it, evolves as an RRSC($m-1$) urn, slowed down by the factor $\lambda_1(m)$. Consequently,}
\begin{equation}
\revv{\delta_1\sim N^{\lambda_1(m-1)/\lambda_1(m)},\quad
\nu(m)=1+\frac{\lambda_1(m)}{\lambda_1(m-1)}}
\label{eq:nu-general}
\end{equation}
\revv{with $\lambda_1(1)=1$ for the RRT. The general prediction \eqref{eq:nu-general} reduces to results stated in Eqs.~(\ref{deg1}), (\ref{deg1-3}), (\ref{deg1-4}) and (\ref{tail-2}), (\ref{tail-3}), (\ref{tail-4}). In particular, the matrix~(\ref{mat:3}) is the transpose of~(\ref{eq:urnM}) for $m=3$, so $\Lambda=\lambda_1(3)$. The exponent $\nu(m)$ is listed in Table~\ref{tab:nu} for $m\le9$; numerically, it decreases towards $\nu(\infty)=2$ as $m\to\infty$.}

\begin{table}
\caption{\revv{Degree exponent of the RRSC($m$) models from Eq.~(\ref{eq:nu-general}).}}
\label{tab:nu}
\revv{%
\begin{tabular*}{\columnwidth}{@{\extracolsep{\fill}}ccccccccc}
\toprule
$m$ & 2 & 3 & 4 & 5 & 6 & 7 & 8 & 9\\
\midrule
$\nu$ & 3.618 & 2.715 & 2.450 & 2.325 & 2.253 & 2.207 & 2.175 & 2.151\\
\bottomrule
\end{tabular*}}
\end{table}

\section{Degree Distribution}
\label{sec:degree}

The computation of the fraction of vertices of degree one is straightforward. We divide the set of all vertices into set of $N_1$ vertices of degree one and set of $N-N_1$ vertices of degree higher than one:
\begin{equation}
\label{sets:S0}
\underbrace{v_i \cdots}_{N_1} \quad \underbrace{v_I \cdots}_{N-N_1}
\end{equation}
with lowcase letters denoting vertices of degree one and capital  letters denoting vertices of degree  higher than one. We divide the set of edges into three subsets:
\begin{equation}
\label{sets:S1-long}
\underbrace{\{v_i,v_j\}\ldots}_{E_1} \quad \underbrace{\{v_k,v_K\}\ldots}_{N_1-2E_1} \quad \underbrace{\{v_I,v_J\}\ldots}_{S_1-N_1+E_1}
\end{equation}
We know $E_1=1$ for $N=2$, see \eqref{K-2}, $E_1=0$ for $N=3$, see \eqref{K3}. Generally for $N\geq 3$, there are no edges between vertices of degree one: $E_1=0$. Therefore, for $N\geq 3$, the division \eqref{sets:S1-long} of the set of edges simplifies to
\begin{equation}
\label{sets:S1}
\underbrace{\{v_k,v_K\}\ldots}_{N_1} \quad \underbrace{\{v_I,v_J\}\ldots}_{S_1-N_1}
\end{equation}

Vertices of degree one do not appear in simplices of dimension $d\geq 2$. When a new vertex $v_{N+1}$ joins a simplex of dimension $d\geq 2$, it acquires degree $d+1$ and the degrees of the vertices from the simplex increase by one and become $\geq d+1$. Therefore, the number of vertices of degree one increases by one if $v_{N+1}$ joins a vertex of the type $v_I$ from \eqref{sets:S0} and decreases by one if $v_{N+1}$ joins an edge of the type $\{v_k,v_K\}$ from \eqref{sets:S1}:
\begin{equation}
\label{N1:eq-2}
N_1\to 
\begin{cases}
N_1+1 &\text{prob} ~~\frac{N-N_1}{S-S_m}\\
N_1-1 &\text{prob} ~~\frac{N_1}{S-S_m}\\
N_1 &\text{prob} ~~1-\frac{N}{S-S_m}
\end{cases}
\end{equation}
Hence, $N_1(N)=\langle N_1\rangle$ satisfies 
\begin{equation}
\label{N1:rec}
N_1(N+1) = N_1(N) + \left\langle \frac{N-2N_1}{S-S_m} \right\rangle
\end{equation}
The linear $\langle N_1\rangle=Nn_1$ growth with $N$ is consistent with \eqref{N1:rec}. One finds $n_1=(1-2n_1)/(1+a_1+\ldots+a_{m-1})$, from 
which the fraction of vertices of degree one is
\begin{equation}
\label{deg:1-m}
n_1 = \frac{1}{3+a_1+\ldots+a_{m-1}} 
\end{equation}
\revv{Using Eq.~\eqref{eq:perron}, we can compactly rewrite Eq.~\eqref{deg:1-m} as $n_1=1/[2+\lambda_1(m)]$.}


For the RRSC(2) model, \eqref{deg:1-m} reduces to the announced result \eqref{deg:1-2} \rev{also matched by numerics (\cref{fig:degree}a)}. Using  \eqref{deg:1-2} and restricting the sums in \eqref{sum:rules} to sums over $\delta\geq 2$ we obtain 
\begin{equation}
\label{sum:rules-2}
\sum_{\delta\geq 2}n_\delta = \frac{15+\sqrt{5}}{22}\,, \qquad  \sum_{\delta\geq 2}(\delta-1) n_\delta = \sqrt{5}
\end{equation}

An extension to vertices of higher degree is perhaps feasible. To determine $n_2$ one should divide the set of all vertices into three sets
\begin{equation}
\label{sets:S0-2}
\underbrace{v_i \cdots}_{N_1} \quad \underbrace{v_\alpha \cdots}_{N_2} \quad \underbrace{v_I \cdots}_{N-N_1-N_2}
\end{equation}
Here, the lowercase Latin letters denote vertices of degree one, lowercase Greek letters denote vertices of degree two, and capital Latin letters denote vertices of degree higher than two. One should further divide the sets of edges and triangles into appropriate subsets. The number of vertices of degree two can change only when vertex $v_{N+1}$ joins a vertex, an edge, or a triangle. It might be possible to derive a stochastic equation for $N_2(N)$ similar to \eqref{N1:rec}, and deduce the fraction $n_2$ of vertices of degree two. The generalization to $n_k$ is in principle straightforward, but the amount of computations rapidly increases with $k$. 

The notion of degree is a natural local characteristic of vertices in graphs: The degree of the vertex $u$ is the number of edges $\{u,v\}$ in the graph adjacent to the vertex $u$. For simplicial complexes, one can consider higher degrees of vertices, accounting for triangles that are adjacent to the vertex, tetrahedra adjacent to the vertex, etc. More generally, if the simplex ${\bf s}$ is the face of $\deg({\bf s},D)$ simplices of dimension $D$ in the simplicial complex $\mathcal{K}$,  we call  $\deg({\bf s},D)$ the degree of ${\bf s}$. 

Among all degrees $\deg({\bf s},D)$ of $d-$dimensional simplices ${\bf s}\in S_d$, we can consider the number of simplices of fixed degree
\begin{equation}
S^{d,D}_\delta = \#[{\bf s}\in \mathcal{K}_d\,| \deg({\bf s},D)=\delta]
\end{equation}

The numbers $S^{d,D}_\delta$ are random in the ensemble of simplicial complexes grown via the RRSC(m) procedure. When $N\gg 1$, these random variables are expected to become self-averaging and scale linearly with $N$. We thus introduce the densities 
\begin{equation}
n^{d,D}_\delta = \lim_{N\to\infty} N^{-1} S^{d,D}_\delta
\end{equation}
There are $\frac{m(m+1)}{2}$ densities $n^{d,D}_\delta$ since dimensions vary in the range $0\leq d<D\leq m$. The density $n^{0,1}_\delta$ is the only one appearing in the RRSC(1). For the RRSC(2), there are three densities: $n^{0,1}_\delta, n^{0,2}_\delta$ and $n^{1,2}_\delta$.

\section{Discussion}
\label{sec:future}

The class of RRSC(m) models admits exact analysis. Understanding models with $m\geq 2$ remains far less complete than understanding the RRSC(1)=RRT model, which has been studied in great depth (see \cite{Pittel94,KR01,KR02-fluct,KR02, Janson05, Drmota,Hofstad,Janson15,Janson19} and references therein). For the RRT model, the numbers of vertices and edges are deterministic: $S_0=N$ and $S_1=N-1$. In Appendix \ref{ap:edges}, we show that $S_1$ and $S_2$ are asymptotically self-averaging for the RRSC(2) model and outline how to extend this treatment and establish the asymptotic self-averaging for $S_1, S_2$, and $S_3$ for the RRSC(3) model. \revv{The asymptotic self-averaging crucial for the analysis presented in Sec.~\ref{sec:linear} follows rigorously for every $m$ from the mapping onto generalized P\'olya urns (Sec.~\ref{sec:urn}), which also reduces the computation of the amplitudes to a linear eigenvalue problem and shows that the fluctuations of $S_d$ are Gaussian with variance $\propto N$ for $m\le5$ and anomalously large for $m\ge6$.} \revv{Numerical simulations for $m\le5$ agree with these predictions.} \revv{Table~\ref{tab:status} summarizes the status of our results.}

\begin{table}
\caption{\revv{Status of the main results.}}
\label{tab:status}
\begin{tabular}{p{0.27\columnwidth}p{0.65\columnwidth}}
\toprule
\raggedright\revv{Exact (all $N$)} & \raggedright\revv{Extreme probabilities (Sec.~\ref{sec:basic}, App.~\ref{ap:edges}); bounds on $S_d$, $\chi=1$; HSC(2): law of $\delta_1$~\cite{DMS01}, fluctuations of $N_2$ (App.~\ref{ap:hom})} \tabularnewline[2pt]
\raggedright\revv{Rigorous, $N\to\infty$} & \raggedright\revv{Amplitudes $a_d(m)$ and self-averaging of $S_d$ for all $m$; Gaussian fluctuations for $m\le5$, growing as $N^{\lambda_2/\lambda_1}$ for $m\ge6$ (Sec.~\ref{sec:urn}, App.~\ref{app:urn})} \tabularnewline[2pt]
\raggedright\revv{Mean-field + numerics} & \raggedright\revv{$\mathrm{Var}[S_1]=N/3$; growth of $\langle\delta_1\rangle$; $n_1$ (Secs.~\ref{sec:primordial}, \ref{sec:degree}, App.~\ref{ap:edges})} \tabularnewline[2pt]
\raggedright\revv{Heuristic + numerics} & \raggedright\revv{$\delta_1\sim N^{\lambda_1(m-1)/\lambda_1(m)}$ with a random prefactor, hence non-self-averaging of $\delta_1$; $\nu(m)=1+\lambda_1(m)/\lambda_1(m-1)$ (Sec.~\ref{sec:nu-general}, App.~\ref{app:urn}, Fig.~\ref{fig:primordial_degree})} \tabularnewline
\bottomrule
\end{tabular}
\end{table}

\revv{Our first computation of the decay exponent characterizing the tail of the degree distribution relied on the scaling law for the growth of the average degree of the primordial vertex and the assumption $\delta_1\sim\delta_{\max}$. The link argument of Appendix~\ref{app:urn} recovers the same exponents for every $m$ without this assumption. Both derivations are heuristic, but such arguments are known to give exact scaling predictions for growing sparse networks,} and they agree well with our simulations, particularly for $m\ge3$. We also computed the fraction of vertices of degree one. Extending this calculation to vertices of degree two is currently a challenge. We also defined higher degree distributions, namely $\frac{m(m+1)}{2}$ densities $n^{d,D}_\delta$ with $0\leq d<D\leq m$. Studying these densities is left for future work. 

\revv{Thanks to the urn mapping, one can determine the amplitudes $a_d(m)$ and the exponent $\nu(m)$ for any $m$; numerically, $\nu(m)$ decreases and approaches $2$ as $m$ grows (we checked up to $m=40$). As $m$ increases, the models change twice. From $m=1$ to $m\ge2$, the main observables change qualitatively: the total number $S_d$ of simplices of dimension $d$ becomes random ($d=1,\ldots,m$), the degree distribution changes from exponential to algebraic, and $\delta_1$ grows as a power of $N$ with a random prefactor, rather than growing as $\ln N$. In the RRT, the link of a vertex grows only through isolated vertices, whereas for $m\ge2$ it grows as an urn of its own (Appendix~\ref{app:urn}). From $m=5$ to $m=6$, only the fluctuations of $S_d$ change, from Gaussian of order $\sqrt{N}$ to order $N^{\lambda_2/\lambda_1}$, because $\lambda_2/\lambda_1$ increases with $m$ and crosses $1/2$; this change appears only asymptotically and is blurred by slow convergence at accessible sizes. The behavior of the amplitudes $a_d(m)$ and of $\lambda_1(m)$ when $m\gg1$ remains an open question. Making the link argument rigorous, e.g., by embedding the link in continuous time as is standard for urns~\cite{AthreyaKarlin68,Janson04}, would establish the growth laws for $\delta_1(m)$ and the exponent $\nu(m)$ for every $m$.} 

The RRSC(m) models provide a natural extension of the RRSC(1), i.e., the RRT model. We now define another class of models that extend the RRT to high-dimensional simplicial complexes and ensure that emerging simplicial complexes are dimensionally uniform. Note that the decomposition \eqref{decomp} can be condensed into
\begin{equation}
\label{decomp-short}
\mathcal{K} = \left\langle \Sigma_1; \ldots; \Sigma_m\right\rangle
\end{equation}
where $\Sigma_j$ contains all $j-$dimensional maximal simplices in $\mathcal{K}$ that are not faces of larger simplices. The RRSC(1) procedure generates simplicial complexes $\left\langle \Sigma_1\right\rangle$ with homogeneous maximal-simplex decomposition and suggests considering the HSC(m) models generating homogeneous simplicial complexes admitting homogeneous maximal-simplex decomposition $\left\langle \Sigma_m\right\rangle$. In Appendix~\ref{ap:hom}, we define the HSC(m) models and show that they are more tractable than the RRSC(m) models. \revv{These models coincide with NGF at $s=0$, $\beta=0$~\cite{Ginestra16-NGF,FIMS22} (and, for $m=2$, with the model of Ref.~\cite{DMS01}); we complement the known degree distribution with exact results on fluctuations and extreme statistics.} The agreement between exact results and predictions of less rigorous methods supports our treatment of the RRSC(m) models. 

Various models of growing deterministic simplicial complexes (DSC) have been investigated, see \cite{Sergey02,Ginestra-Ziff,Ginestra-Sergey20,Ginestra-reitz,Sergey25} and references therein. One particularly close study \cite{Sergey25} explores the family of constrained DSC(m) models in which the dimension of simplices does not exceed $m$. It would be interesting to extend computations \cite{Sergey25} of the Hodge Laplacian spectra, upper degree distributions, and other properties of deterministic simplicial complexes generated by the DSC(m) models to random simplicial complexes generated by the RRSC(m) models. 

\bigskip\noindent
We are grateful to Ginestra Bianconi and Sergey Dorogovtsev for discussions.
\rev{M.L. was a postdoctoral researcher at the Fonds de la Recherche Scientifique-FNRS.}

\appendix

\section{Total numbers of simplices}
\label{ap:edges}

Our analytical work was based on the asymptotic self-averaging of the total number of simplices of various dimensions. In this Appendix, we justify this assumption for the RRSC(2) model and outline how to extend the arguments to the RRSC(3) model. \revv{A general proof for all $m$ follows from the urn mapping of Sec.~\ref{sec:urn}; here we present an elementary derivation that also yields the variance explicitly.} 

Strictly speaking, the analysis presented below is complete only for the RRSC(2) model, in which we can focus on a single random variable, e.g., the number of edges $S_1$. Specializing the Euler characteristic \eqref{Euler:def} to $m=2$, we find the relation
\begin{equation}
S_1 - S_2 = N-1
\end{equation}
showing that the random number of triangles $S_2$ is enslaved to the number of edges. 

To study the evolution of the number of edges, we note that adding a new vertex, $N\to N+1$, leads to the increase of the number of edges 
\begin{equation}
\label{S1:rec}
S_1\to 
\begin{cases}
S_1+1          & \text{prob} \quad \frac{N}{N+S_1}\\
S_1+2           &\text{prob} \quad \frac{S_1}{N+S_1}
\end{cases}
\end{equation}
Despite the neat form, the stochastic growth law \eqref{S1:rec} is non-linear and hence challenging. 

The probability distribution $P(S_1,N)$ satisfies the difference equation
\begin{eqnarray}
\label{PS1}
P(S_1,N+1) &=& \frac{N}{N+S_1-1}\,P(S_1-1,N)\nonumber \\
&+&  \frac{S_1-2}{N+S_1-2}\,P(S_1-2,N)
\end{eqnarray}
immediately following from \eqref{S1:rec}. The probabilities of generating the minimal number of edges, $S_1=N-1$, and the maximal number of edges, $S_1=2N-3$, easily follow from \eqref{PS1}. Note that $P(N-1, N)$ coincides with the probability $T_N$ of generating a tree given by Eq.~\eqref{TN-sol}. The probability of generating simplicial complexes with the maximal number of edges, i.e., to build an open book, satisfies the recurrence
\begin{equation}
\label{PS1-max}
P(2N-1,N+1)=\frac{2N-3}{3N-3}\,P(2N-3,N)
\end{equation}
following from \eqref{PS1}. Iterating \eqref{PS1-max} yields 
\begin{align}
\label{max:S1}
P(2N-3,N)= \left(\frac{2}{3}\right)^{N-2} \frac{\Gamma(N-\frac{3}{2})}{\sqrt{\pi}\,\Gamma(N-1)}
\end{align}

\begin{figure*}[htb]
	\centering
	\includegraphics[width=1\linewidth]{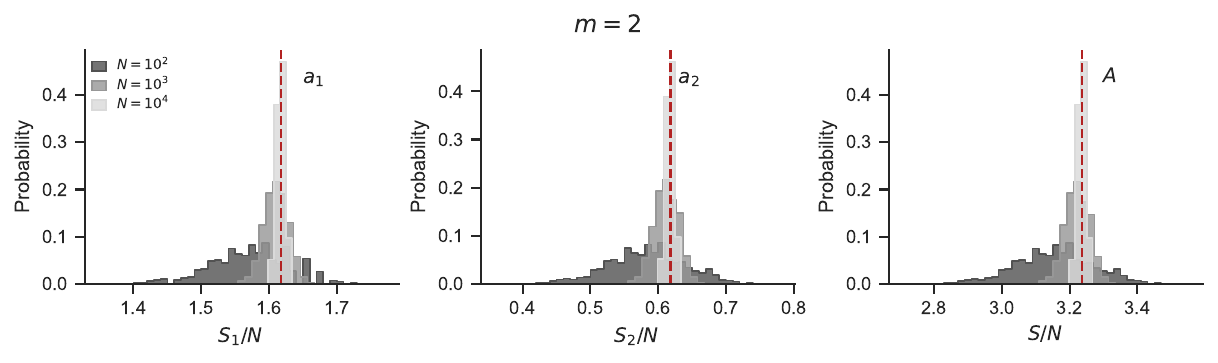}
	\caption{\revv{Concentration and self-averaging of the $S_d$ for RRSC(2).}
		Distributions of all $S_d/N$ and $S/N$ at $N=10^2$, $10^3$, and $10^4$ concentrate sharply around  corresponding $a_d$ and $A$ (red dashed) from \cref{a1:2} as $N$ grows.}
	\label{fig:Sd_m2}
\end{figure*}

\begin{figure*}[htb]
	\centering
	\includegraphics[width=1\linewidth]{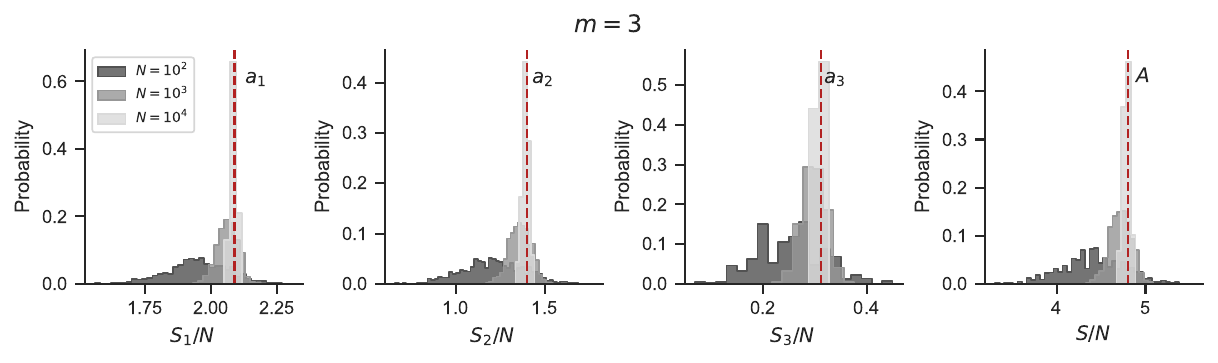}
	\caption{\revv{Concentration and self-averaging of the $S_d$ for RRSC(3).}
		Distributions of all $S_d/N$ and $S/N$ at $N=10^2$, $10^3$, and $10^4$ concentrate sharply around  corresponding $a_d$ and $A$ (red dashed) from \cref{aA:3} as $N$ grows.}
	\label{fig:Sd_m3}
\end{figure*}

\begin{figure*}[htb]
	\centering
	\includegraphics[width=1\linewidth]{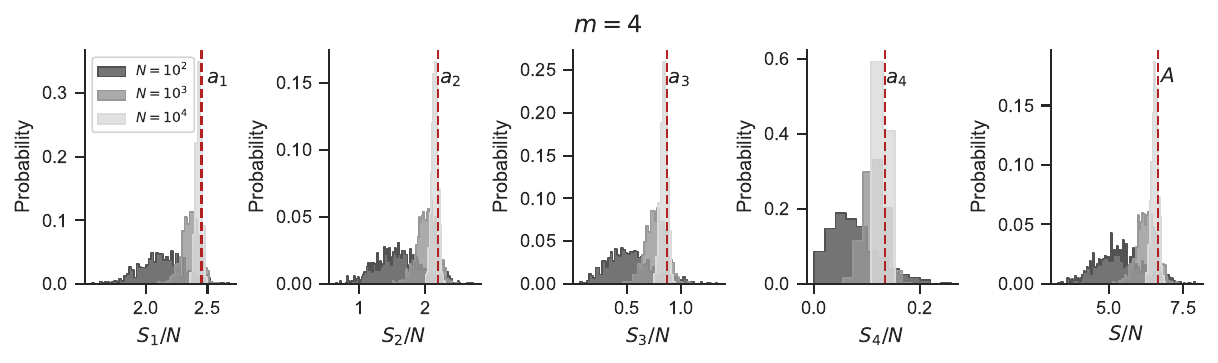}
	\caption{\revv{Concentration and self-averaging of the $S_d$ for RRSC(4).}
		Distributions of all $S_d/N$ and $S/N$ at $N=10^2$, $10^3$, and $10^4$ concentrate sharply around  corresponding $a_d$ and $A$ (red dashed) from \cref{aA:4} as $N$ grows.}
	\label{fig:Sd_m4}
\end{figure*}

One would like to solve \eqref{PS1} subject to the initial condition $P(S_1,2)=\delta_{S_1,1}$, but it is impossible because the amplitudes on the right-hand side of \eqref{PS1} are non-linear functions of $S_1$. Fortunately, one can determine the asymptotic behavior of the variance $\text{Var}[S_1]$ as we now demonstrate. Multiplying \eqref{PS1} by $S_1$ and summing over all $S_1$ we obtain 
\begin{equation}
\label{S1:eq-1}
\langle S_1(N+1)\rangle = \langle S_1\rangle + 2-Nv
\end{equation}
where 
\begin{equation}
\label{v:def}
v=\left\langle \frac{1}{N+S_1}\right\rangle
\end{equation}
and we shortly write $S_1$ instead of $S_1(N)$; we use the complete notation only when the number of vertices is $N+1$. Multiplying \eqref{PS1} by $S_1^2$ and summing over all $S_1$ we deduce 
\begin{equation}
\label{S1:eq-2}
\langle S_1^2(N+1)\rangle = \langle S_1^2\rangle -2N + 4\langle S_1\rangle + 4+ 2N^2 v - 3 Nv 
\end{equation}
Subtracting the square of \eqref{S1:eq-1} from \eqref{S1:eq-2} gives
\begin{eqnarray}
\label{S1:eq-var}
\text{Var}[S_1(N+1)] &=& \text{Var}[S_1] + Nv - (N v)^2 \nonumber\\
&+&2N\left[Nv + \langle S_1\rangle v-1\right]
\end{eqnarray}

None of the equations \eqref{S1:eq-1}, \eqref{S1:eq-2}, and \eqref{S1:eq-var} is a recurrence due to the presence of $v$. However, the progress is feasible when $N\gg 1$. In this situation, we can decompose $S_1$ into a deterministic part, the average that scales $N$, and a random part that scale as $\sqrt{N}$:
\begin{equation}
\label{sigma:def}
S_1 = \langle S_1\rangle + \sigma \sqrt{N}
\end{equation}
The random quantity $\sigma$ becomes independent of $N$ when $N\to\infty$. We now insert \eqref{sigma:def} into the term in the square bracket in \eqref{S1:eq-var} and recast it into
\begin{eqnarray*}
Nv + \langle S_1\rangle v-1 &=& \left\langle \frac{N+\langle S_1\rangle}{N+\langle S_1\rangle+\sigma \sqrt{N}}\right\rangle-1 \\
&=& -\left\langle \frac{\sigma \sqrt{N}}{N+\langle S_1\rangle+\sigma \sqrt{N}}\right\rangle 
\end{eqnarray*}
Expanding we obtain 
\begin{equation}
\frac{N}{(N+\langle S_1\rangle)^2}\,\langle \sigma^2\rangle - \frac{N^{3/2}}{(N+\langle S_1\rangle)^3}\,\langle \sigma^3\rangle+\ldots
\end{equation}
where we have used $\langle \sigma\rangle = 0$ following from \eqref{sigma:def}. Recalling the asymptotic $\langle S_1\rangle\simeq a_1N$ with $a_1=\frac{1+\sqrt{5}}{2}$ we establish the leading behaviors $Nv \simeq (1+a_1)^{-1}$ and 
\begin{equation*}
\begin{split}
2N\left[Nv + \langle S_1\rangle v-1\right]  &\simeq \frac{2N^2}{(N+\langle S_1\rangle)^2}\,\langle \sigma^2\rangle\\
& = \frac{2}{(1+a_1)^2}\,\langle \sigma^2\rangle 
\end{split}
\end{equation*}
Plugging these asymptotic behaviors into \eqref{S1:eq-var} and using $\text{Var}[S_1] = N\langle \sigma^2\rangle$ following from \eqref{sigma:def} we obtain 
\begin{equation*}
\langle \sigma^2\rangle = (1+a_1)^{-1}-(1+a_1)^{-2}+ \frac{2}{(1+a_1)^2}\,\langle \sigma^2\rangle
\end{equation*}
which we solve and find $\langle \sigma^2\rangle=\frac{1}{3}$ as announced in \eqref{S1:var}. 

Equations \eqref{S1:eq-1}, \eqref{S1:eq-2}, and \eqref{S1:eq-var} are exact, and one can utilize them for computing sub-leading corrections. We illustrate how to deal with the average. We write
\begin{equation}
\label{eps:def}
\langle S_1\rangle = a_1N + \epsilon(N)
\end{equation}
where $\epsilon(N)$ is an unknown sub-leading correction. Inserting \eqref{sigma:def} and \eqref{eps:def} into \eqref{v:def} 
we obtain
\begin{eqnarray}
\label{Nv}
Nv &=& \left\langle \frac{N}{N+a_1N + \epsilon(N) + \sigma \sqrt{N}}\right\rangle \nonumber \\
& = & - \frac{\epsilon}{(1+a_1)^2 N} + \frac{\langle \sigma^2\rangle}{(1+a_1)^3 N} + \ldots
\end{eqnarray}
Inserting this result into \eqref{S1:eq-1} gives
\begin{equation}
\label{eps:eq}
\frac{d\epsilon}{dN} =\frac{\epsilon}{(1+a_1)^2 N}
\end{equation}
The replacement of the difference $\epsilon(N+1)-\epsilon(N)$ by the derivative $\frac{d\epsilon}{dN}$ is justified in the $N\to\infty$ limit. Solving \eqref{eps:eq} gives
\begin{equation}
\label{eps:sol}
\epsilon \sim N^{1/(1+a_1)^2}\,, \quad \frac{1}{(1+a_1)^2} \approx 0.145898
\end{equation}
The exponent is small but positive, so the first term from \eqref{Nv} dominates the second, explaining why we kept only the first term in \eqref{eps:eq}. \revv{The exponent $(1+a_1)^{-2}=(7-3\sqrt5)/2$ equals $\lambda_2/\lambda_1$ for the matrix~\eqref{eq:urnM} with $m=2$.} 

A straightforward but tedious analysis of Eq.~\eqref{PS1} shows that the distribution $P(S_1,N)$ is asymptotically Gaussian, and therefore fully characterized by the average and the variance
\begin{equation}
\label{PSN:Gauss}
P(S_1,N)= \sqrt{\frac{3}{2\pi\,N}}\,\exp\!\left[-\frac{3(S_1-a_1N)^2}{2 N}\right]
\end{equation}

The extension to the RRSC(m) models is, in principle, straightforward but requires laborious calculations. For the RRSC(3) model, we can choose $S_1$ and $S_2$ as basic random variables. The number of tetrahedrons $S_3$ can be expressed via $S_1$ and $S_2$:
\begin{equation}
S_3 = N - 1 - S_1 + S_2
\end{equation}
Adding a new vertex, $N\to N+1$, leads to 
\begin{equation}
\label{S12:recur}
(S_1,S_2)\to 
\begin{cases}
(S_1+1,S_2)               & \text{prob} \quad \frac{N}{N+S_1+S_2}\\
(S_1+2,S_2+1)           &\text{prob} \quad \frac{S_1}{N+S_1+S_2} \\
(S_1+3,S_2+3)           &\text{prob} \quad \frac{S_2}{N+S_1+S_2}
\end{cases}
\end{equation}
The growth law \eqref{S12:recur} implies that the probability distribution $P(S_1,S_2;N)$ satisfies
\begin{eqnarray}
\label{PS12}
P(S_1,S_2;N+1) &=& \frac{S_2-3}{M-6}\,P(S_1-3,S_2-3;N) \nonumber \\
                            &+&  \frac{S_1-2}{M-3}\,P(S_1-2,S_2-1;N) \nonumber \\
                            &+& \frac{N}{M-1}\,P(S_1-1,S_2;N)
\end{eqnarray}
where $M=N+S_1+S_2$. One can again compute the extreme probabilities. We know $P(N-1, 0; N)=T_N$, with $T_N$ being the probability to generate a tree, Eq.~\eqref{TN-sol}. The probability to generate the maximal numbers of edges and triangles
\begin{align*}
P(3N-6,3N-8; N)= \left(\frac{3}{7}\right)^{N-3} \frac{\Gamma(N-\frac{8}{3})}{3\Gamma(\frac{1}{3})\,\Gamma(N-2)}
\end{align*}
is found by iterating \eqref{PS12}.

Decomposing $S_1$ and $S_2$ into deterministic parts varying as $N$ and random parts scaling as $\sqrt{N}$,
\begin{equation}
\label{sigma:12}
S_1 = \langle S_1\rangle + \sigma_1 \sqrt{N}\,, \quad S_2 = \langle S_2\rangle + \sigma_2 \sqrt{N}\,, 
\end{equation}
and substituting these expressions into equations for the variances $\text{Var}[S_j] = N\langle \sigma_j^2\rangle$ with $j=1,2$ and the covariance $\langle\!\langle S_1S_2\rangle\!\rangle=\langle S_1S_2\rangle-\langle S_1\rangle \langle S_2\rangle=N\langle \sigma_1 \sigma_2\rangle$, one arrives at a system of linear equations for $\langle\sigma_1^2\rangle$, $\langle\sigma_2^2\rangle$, and $\langle \sigma_1 \sigma_2\rangle$, which are solved to yield
\begin{equation}
\langle \sigma_1^2\rangle = \lambda_1, \quad  \langle \sigma_2^2\rangle = \lambda_2 , \quad \langle \sigma_1 \sigma_2\rangle = \mu 
\end{equation}
We omit cumbersome explicit expressions of the variances $\lambda_1, \lambda_2$, and the covariance $\mu$. (Recall, that even the amplitudes $a_1$ and $a_2$ are not solvable in quadratures for the RRSC(3) model.)

The probability distribution becomes stationary as a function of $\sigma_1$ and $\sigma_2$, i.e., $P(S_1,S_1; N)=\Pi(\sigma_1,\sigma_2)$, and acquires a Gaussian form
\begin{equation}
\Pi=\frac{1}{2\pi\sqrt{\lambda_1 \lambda_2- \mu^2}}\,\exp\!\left\{-\frac{F(\sigma_1,\sigma_2)}{2(\lambda_1 \lambda_2- \mu^2)}\right\}
\end{equation}
with 
\begin{equation}
F(\sigma_1,\sigma_2) = \lambda_2 \sigma_1^2+\lambda_1 \sigma_2^2-2\mu  \sigma_1  \sigma_2
\end{equation}

\revv{Figures~\ref{fig:Sd_m2}, \ref{fig:Sd_m3}, and~\ref{fig:Sd_m4} corroborate the theoretical results of Sec.~\ref{sec:linear}: for $m=2$, $3$, and $4$, the distributions of all $S_d/N$ and of $S/N$ concentrate around the amplitudes $a_d$ and $A$ given by Eqs.~(\ref{a1:2}), (\ref{aA:3}), and~(\ref{aA:4}) as $N$ grows.}

\section{Amplitudes $a_d$}
\label{ap:10}

Extending the analysis presented in Sec.~\ref{sec:linear} to the RRSC(m) models with $m\geq 5$ is straightforward. In Table~\ref{Table:10}, we collect the amplitudes $a_d(m)$ for the RRSC(m) models with $m < 10$. The amplitudes can be determined with any desirable precision. The sum rule
\begin{equation}
1+\sum_{d=1}^m (-1)^d a_d(m) = 0
\end{equation}
follows from the identity \eqref{Euler:def} for the Euler characteristic and the linear growth of $S_d$. 


\begin{table}[h!]
\centering
\renewcommand{\arraystretch}{1.5}{
\begin{tabular}{| c | c | c | c | c | c | c | c | c | c |}
\hline
$m$         &  $a_1$     &  $a_2$      &  $a_3$      &  $a_4$   & $a_5$   & $a_6$  & $a_7$  & $a_8$ & $a_9$\\ 
\hline
1              & 1             & 0                  & 0              & 0            & 0          & 0          & 0          & 0          & 0  \\ 
\hline
2              & 1.618      & 0.618           & 0              & 0            & 0          & 0          & 0          & 0          & 0    \\ 
\hline
3              & 2.089      & 1.401           & 0.312        & 0            & 0          & 0         & 0          & 0          & 0     \\ 
\hline
4              & 2.450      & 2.188           & 0.871        & 0.133      & 0          & 0         & 0          & 0          & 0     \\ 
\hline
5              & 2.731      & 2.908           & 1.556        & 0.429      & 0.050    & 0         & 0          & 0          & 0    \\ 
\hline
6              & 2.957      & 3.551           & 2.280        & 0.846      & 0.176   & 0.016   & 0         & 0          & 0   \\
\hline
7              & 3.146      & 4.126           & 3.002       & 1.339      & 0.374    & 0.061   & 0.004   & 0          & 0   \\
\hline
8              & 3.307      & 4.645           & 3.708       & 1.878      & 0.630     & 0.140   & 0.019   & 0.001   & 0   \\
\hline
9              & 3.447      & 5.12           & 4.394        & 2.447      & 0.934     & 0.248   & 0.045   & 0.005     &  0.0003\\
\hline
\end{tabular}
}
\caption{The amplitudes $a_d(m)$ for the RRSC(m) models with $m<10$. The amplitudes with $d\geq 10$ are trivial, $a_d(m)=0$. 
The $10^{-3}$ precision threshold is used for data shown in the Table, apart from $a_9(9)\approx 2.97\cdot 10^{-4}$ which is smaller than $10^{-3}$.}
\label{Table:10}
\end{table}

\begin{figure}[hb]
	\centering
	\includegraphics[width=1\linewidth]{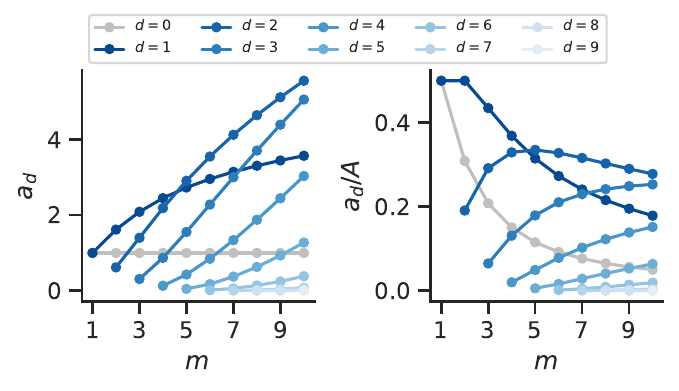}
	\caption{The amplitudes $a_d(m)$ for the RRSC(m) models with $m<10$ from \cref{Table:10}.}
	\label{fig:ad}
\end{figure}

The values of the amplitudes listed in Table~\ref{Table:10} \rev{and shown in \cref{fig:ad}} lead to a few `experimental' observations:
\begin{enumerate}
\item By definition of the RRSC(m) models, $a_d(m)=0$ when $m<d$. For $m\geq d$ and fixed $d$,  the amplitudes $a_d(m)$ monotonically increase with $m$. 
\item An algebraic growth, $a_d(m)\sim m^{\alpha_d}$, seem plausible when $m\gg 1$. The data presented in Table~\ref{Table:10} is too small for estimating the exponents; we mention $\alpha_1\approx 0.5$ and $\alpha_2\approx 0.75$ as the simplest guesses. 
\item The amplitudes $a_m(m)$ corresponding to simplices of the largest dimension quickly decrease with $m$; an exponential decay seems plausible.  
\item The amplitudes $a_{m-k}(m)$ with fixed $k$ quickly, probably exponentially, decrease with $m$. 
\end{enumerate}

\section{Mapping onto a generalized P\'olya urn}\label{app:urn}

\revv{We now explain the mapping of the RRSC($m$) process onto a generalized P\'olya urn and derive the results summarized in Secs.~\ref{sec:urn} and~\ref{sec:nu-general}.}

\revv{In generalized P\'olya urns, an urn contains balls of several types. At each step, a ball is drawn uniformly at random and then returned, together with new balls; the number of added balls of each type depends only on the type of the drawn ball. This rule is encoded in a replacement matrix whose entry $(d,k)$ is the number of type-$k$ balls added when a type-$d$ ball is drawn. The classical P\'olya urn, with two colors and one ball of the drawn color added at each step, is the simplest example. Generalized urns are a standard tool for proving limit theorems for random recursive structures; see Refs.~\cite{Mahmoud08,FlajoletDumasPuyhaubert06} for introductions and Refs.~\cite{AthreyaKarlin68,Janson04} for the general theory. The long-time behavior of an urn is governed by the eigenvalues and eigenvectors of its replacement matrix.}

\revv{To map our system, we consider each simplex of dimension $d<m$ as a ball of type $d$. Choosing a simplex uniformly among the simplices of dimension smaller than $m$ is the same as drawing a ball uniformly, and the chosen simplex stays in the complex, just as the drawn ball is put back. Let $\sigma$ be the chosen $d$-simplex and $v$ the new vertex. The faces of the new simplex $\{v\}\cup\sigma$ can be divided into two groups. First, the faces that do not contain $v$ are faces of $\sigma$, and are therefore already present. Second, the faces that contain $v$ take the form $\{v\}\cup\tau$ with $\tau\subseteq\sigma$, and are all new, since $v$ did not exist before. A $k$-dimensional face of this kind corresponds to a choice of $k$ vertices of $\sigma$, so drawing a $d$-simplex adds exactly $\binom{d+1}{k}$ simplices of dimension $k$, for $k=0,\ldots,d+1$. For example, choosing a triangle adds one vertex, three edges, three triangles, and one tetrahedron. Simplices of dimension $m$ are created but never chosen, so they are not balls and are counted separately. This gives the replacement matrix~(\ref{eq:urnM}), with rows labelling the type drawn and columns the type added. For $m=3$, for instance, its rows are $(1,1,0)$, $(1,2,1)$, and $(1,3,3)$. The replacements are deterministic; the only randomness lies in which type is drawn.}

\revv{The system can now be reduced to an eigenvector problem. The matrix $M$ is non-negative. It is also irreducible, i.e., every type eventually produces every other type, since drawing a $d$-simplex always creates a vertex ($M_{d0}=1$) and a $(d+1)$-simplex ($M_{d,d+1}=1$ for $d<m-1$). By the Perron--Frobenius theorem \cite[Chap.~8]{HornJohnson13}, the largest eigenvalue $\lambda_1$ of $M$ is real and non-degenerate; its left eigenvector, with components $w_d$ satisfying $\sum_dw_dM_{dk}=\lambda_1w_k$, can be chosen positive, and every other eigenvalue has real part smaller than $\lambda_1$. For urns with finitely many types, equal drawing weights, an irreducible replacement matrix, and replacements of finite variance, all of which hold here, the composition of the urn grows linearly with a deterministic limit~\cite{AthreyaKarlin68}, \cite[Theorem~3.21]{Janson04}:}
\begin{equation}
\revv{\lim_{N\to\infty}\frac{S_d(N)}{N}=\lambda_1 w_d}
\label{eq:LLN}
\end{equation}
\revv{almost surely, with $\sum_{0\leq d<m}w_d=1$. The total number of balls, $\sum_{d<m}S_d$, thus grows as $\lambda_1N$, and the fraction $S_d/\sum_{d'<m}S_{d'}$ of balls of type $d$ tends to $w_d$.}

\revv{We can now express this in terms of our amplitudes $a_d$, which were defined in Eq.~(\ref{Sd:av}) through the growth of the average, $\langle S_d\rangle\simeq a_dN$. Each step adds at most $2^m$ simplices, so $S_d/N\le 2^m$, and the almost sure convergence~(\ref{eq:LLN}) implies, by dominated convergence, $\langle S_d\rangle/N\to\lambda_1w_d$. Hence the limit in Eq.~(\ref{Sd:av}) exists and $a_d=\lambda_1w_d$, so the amplitudes satisfy $\lambda_1a_k=\sum_{d}a_d\binom{d+1}{k}$. Since each step adds exactly one vertex, $a_0=1$, and the $k=0$ component gives $\lambda_1=1+a_1+\cdots+a_{m-1}$, which is Eq.~(\ref{eq:perron}). With this value of $\lambda_1$, the components $k=1,\ldots,m-1$ are precisely Eqs.~(\ref{ak-m}). Moreover, each draw of an $(m-1)$-simplex creates one $m$-simplex, which happens with probability $a_{m-1}/\lambda_1$; this recovers Eq.~(\ref{am-m}). Since the Perron eigenvector is the only eigenvector with positive components, the positive solution of Eqs.~(\ref{ak-m}) is unique. For $m=2$, $M=\bigl(\begin{smallmatrix}1&1\\1&2\end{smallmatrix}\bigr)$ has $\lambda_1=(3+\sqrt5)/2$, giving $a_1=\lambda_1-1=(1+\sqrt5)/2$, as in Eq.~(\ref{a1:2}). For larger $m$, diagonalizing $M$ reproduces the amplitudes to all digits shown Appendix~\ref{ap:10}.}

\revv{The same bound also yields self-averaging: by dominated convergence, Eq.~(\ref{eq:LLN}) implies convergence in mean square, $\langle(S_d/N-a_d)^2\rangle\to0$; hence $\mathrm{Var}[S_d]/\langle S_d\rangle^2\to0$ for every $d\le m$ and every $m$, which justifies the assumption made in Secs.~\ref{sec:linear}--\ref{sec:degree}.}

\begin{figure}[t]
	\centering
	\includegraphics[width=1\linewidth]{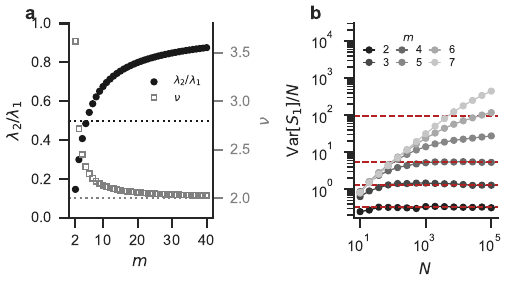}
	\caption{\revv{Generalized P\'olya urn description of the RRSC($m$) models. (a) Ratio $\lambda_2/\lambda_1$ of the two largest eigenvalues of the matrix~(\ref{eq:urnM}) and degree exponent $\nu$ from Eq.~(\ref{eq:nu-general}) as functions of $m$; dotted lines mark $1/2$ and $2$. (b) $\mathrm{Var}[S_1]/N$ for $m=2,\ldots,7$, from simulations of the urn with $R=1000$ independent realizations. Dashed lines: asymptotic values for $m\le5$, where $\lambda_2<\lambda_1/2$, from the linear-noise approximation ($1/3$ for $m=2$, Eq.~(\ref{S1:var})). For $m\ge6$, the variance grows faster than $N$.}}
	\label{fig:urn}
\end{figure}

\revv{For such urns, the fluctuations around $a_dN$ are controlled by the ratio $r=\lambda_2/\lambda_1$~\cite{Janson04,Mahmoud08}. Heuristically, the projection $y_j$ of the composition onto the right eigenvector of $M$ with eigenvalue $\lambda_j$ has an average drift $(\lambda_j/\lambda_1)\,y_j/N$, so a deviation present at step $n$ is amplified by $(N/n)^{\lambda_j/\lambda_1}$ by step $N$; the largest such exponent, for $j\ge2$, is $r$. This sets the threshold $r=1/2$. Indeed, each draw also adds to $y_j$ a random kick with zero mean and bounded variance, uncorrelated with earlier kicks. The deviation at step $N$ is the sum of these kicks, each amplified by $(N/n)^{r}$, so its variance scales as $N^{2r}\sum_{n\le N}n^{-2r}$. For $r<1/2$, the sum scales as $N^{1-2r}$; therefore, the variance grows linearly in $N$, and $(S_d-a_dN)/\sqrt{N}$ converges to a Gaussian~\cite[Theorem~3.22]{Janson04}. For $r>1/2$, the variance is dominated by the earliest kicks, and generically $(S_d-a_dN)/N^{r}$ converges to a random limit set by the early history~\cite[Theorem~3.24]{Janson04}. The same drift makes the average approach $a_dN$ with a correction of order $N^{r}$, set by the initial condition. For $m=2$, $r=(7-3\sqrt5)/2=(1+a_1)^{-2}$, which is precisely the exponent found in Eq.~(\ref{eps:sol}). For our matrices, all eigenvalues are real, and $r$ increases with $m$ (we verified this numerically up to $m=40$), crossing $1/2$ between $m=5$ and $m=6$ (Table~\ref{tab:ratio} and Fig.~\ref{fig:urn}a). The threshold is sharp only asymptotically. Figure~\ref{fig:urn}b shows that $\mathrm{Var}[S_1]/N$ saturates at its predicted asymptotic value for $m\le4$, but keeps growing over the accessible range of $N$ for $m\ge6$, with a slow approach to the asymptotic growth $N^{2r}$. For $m=5$, where $r$ is close to $1/2$, the variance converges to its asymptotic value only as $N^{2r-1}\approx N^{-0.03}$, far beyond the simulated sizes.}

\revv{The urn picture also yields the degree of a vertex for every $m$, through its link, and explains the results of Sec.~\ref{sec:primordial}. The link of a vertex $x$ is the set of simplices $\tau$ such that $x\notin\tau$ and $\tau\cup\{x\}$ belongs to the complex; it is itself a simplicial complex~\cite{Munkres84}. Each $d$-simplex $\tau\cup\{x\}$ containing $x$ thus corresponds to the $(d-1)$-simplex $\tau$ of the link, the empty face $\emptyset$ corresponding to $x$ itself. Its vertices are the neighbours of $x$, so the degree of $x$ is the number of vertices of its link. Any set of pairwise adjacent vertices of an RRSC is a simplex, since its youngest vertex joined a simplex containing all the others. The link of $x$ is therefore its ego network with $x$ removed and all cliques filled. In Fig.~\ref{fig:urnscheme}a, for instance, the link of $c$ consists of the isolated vertex $b$ and the edge $\{d,e\}$, and $c$ has degree $3$.}

\revv{To follow the degree of $x$, we therefore follow how its link grows. When a new vertex $v$ joins a simplex $\sigma$, the link of $x$ changes only if $x\in\sigma$. Writing $\sigma=\{x\}\cup\tau$, the link then acquires $\{v\}\cup\tau'$ for all $\tau'\subseteq\tau$: inside the link, $v$ joins the simplex $\tau$, of dimension $<m-1$, which is uniformly distributed among such elements of the link since the global choice is uniform. This is the RRSC($m-1$) rule, applied to a link that need not be connected, except for $\tau=\emptyset$ (i.e., $\sigma=\{x\}$), where $v$ is added as an isolated vertex. In Fig.~\ref{fig:urnscheme}c, the link of $b$ acquires $\{v\}$ and $\{v,a\}$, i.e., $v$ joins the vertex $a$ of the link, while the link of $c$ is unchanged. Through this correspondence, the elements of the link of dimension $<m-1$, including $\emptyset$, are the balls of the urn that contain $x$, and each is drawn with probability $\simeq1/(\lambda_1(m)N)$. The link of $x$ hence evolves as an RRSC($m-1$) urn that is updated only when one of its balls is drawn. Each move, i.e., each step at which the chosen simplex contains $x$, adds one vertex to the link, so that $\delta=\delta_0+n$, where $\delta_0$ is the degree of $x$ at birth and $n$ the number of such moves. After $n$ moves, the link holds $\simeq\lambda_1(m-1)\,n$ balls, so $dn/dN\simeq[\lambda_1(m-1)/\lambda_1(m)]\,n/N$. Hence $n$, and thus $\delta$, grows as $N^{\lambda_1(m-1)/\lambda_1(m)}$. The isolated vertices added through $\tau=\emptyset$ contribute only $O(\ln N)$ and do not affect this exponent. This gives the first relation in Eq.~(\ref{eq:nu-general}). The argument is heuristic, at the same mean-field level as Sec.~\ref{sec:primordial}, but it holds for every $m$ and explains the results found there for $m\le4$. The quantities $\delta_1$, $\Delta_1$, and $\tau_1$ are the numbers of vertices, edges, and triangles of the link of the primordial vertex. Equations~(\ref{delta-eq}), (\ref{delta-Delta}), and (\ref{delta-Delta-tau}) are thus the mean evolution equations of an RRSC($m-1$) urn, which is why $\Lambda=\lambda_1(3)$ and why the ratios~(\ref{d-D:ratio}) and~(\ref{two-ratios}) are the amplitudes $a_1(2)$, and $a_1(3)$ and $a_2(3)$, respectively. Since the link starts with a few elements, its first moves are random and set the prefactor of this growth law, which remains random, consistent with the lack of self-averaging of $\delta_1$ in Fig.~\ref{fig:primordial_degree}. The exponent $\nu$ then follows without assuming $\delta_1\sim\delta_{\max}$: a vertex born at time $N_0=uN$, with $u$ uniform in $(0,1)$, has degree $\sim u^{-\lambda_1(m-1)/\lambda_1(m)}$ at time $N$, so that $\mathrm{Prob}[\delta>k]\sim k^{-\lambda_1(m)/\lambda_1(m-1)}$, which yields the second relation in Eq.~(\ref{eq:nu-general}). Numerically, $\nu$ decreases monotonically towards $2$ as $m$ grows (Fig.~\ref{fig:urn}a). Finally, using Eq.~(\ref{eq:perron}), the fraction~(\ref{deg:1-m}) of vertices of degree one can be written as $n_1=1/[2+\lambda_1(m)]$.}

\section{Homogeneous simplicial complexes}
\label{ap:hom}

The natural recursive process generating homogeneous simplicial complexes of dimension $m$ begins with an $m-$dimensional simplex and proceeds by adding, on each step, a new vertex $v$ and a new $m-$dimensional simplex $\{v,{\bf s}\}$ where ${\bf s}$ is an $(m-1)-$dimensional existing simplex chosen uniformly among such simplices. \revv{The HSC($m$) model coincides with the network geometry with flavor (NGF)~\cite{Ginestra16-NGF} with flavor $s=0$ and inverse temperature $\beta=0$ in dimension $d=m$, equivalently with Model~A of Ref.~\cite{FIMS22} with constant fitness. For $m=2$, its one-dimensional skeleton is the growing network introduced in Ref.~\cite{DMS01}, in which each new vertex attaches to both ends of a uniformly chosen edge. The limiting degree distribution is $P(K)=12/[K(K+1)(K+2)]$ for $m=2$~\cite{DMS01}. The tail of the degree distribution is known for all $m$: $P(K)\sim K^{-\nu}$ with $\nu=2+1/(m-1)$~\cite{Ginestra16-NGF,FIMS22}. Below, we recover these results with our methods and focus on quantities that, to our knowledge, have not been studied for these models: the fluctuations of the number of vertices of minimal degree and extreme-value statistics.} 

The class of HSC(m) models is more tractable than the class of RRSC(m) models. When $m=1$, the HSC(1) and RRSC(1) are identical to the RRT model. For the HSC(m) models with $m\geq 2$, the total number of simplices of each dimension remains deterministic:
\begin{equation}
\label{Sd:N}
S_d=\binom{m}{d}N-d\binom{m+1}{d+1}
\end{equation}
for $d=0,1,\ldots,m$. This feature simplifies analysis: We do not need to worry about asymptotic self-averaging and extensivity. The total number of simplices is also a deterministic quantity:
\begin{equation}
S = 2^m N - (m-1) 2^m -1
\end{equation}

As an illustration, let us consider the HSC(2). The process begins with a single triangle and 
\begin{equation}
S_0=N, \quad S_1=2N-3, \quad S_2=N-2
\end{equation}
throughout the evolution. \revv{The emerging simplicial complex is a pure two-dimensional complex made of triangles glued along edges; the open book (Fig.~\ref{fig:book_illust}) is one particular realization.}

The primordial triangle is unique, while any of its three vertices can be called primordial. We take one of them and follow the evolution of the degree $\delta_1$ of the primordial vertex. Adding a vertex, $N\to N+1$, leads to the possible increase of the degree $\delta_1$:
\begin{equation}
\label{delta-HSC-2}
\delta_1\to 
\begin{cases}
\delta_1+1 & \text{prob} \quad \frac{\delta_1}{2N-3}\\
\delta_1 & \text{prob} \quad 1-\frac{\delta_1}{2N-3}
\end{cases}
\end{equation}
Using \eqref{delta-HSC-2} we deduce the recurrence relation
\begin{equation}
\label{av-HSC}
\delta_1(N+1)=\left[1+\frac{1}{2N-3}\right]\delta_1(N)
\end{equation}
for the average $\delta_1(N)=\langle \delta_1\rangle$. Starting with $\delta_1(3)=2$ and iterating \eqref{av-HSC} we obtain
\begin{equation}
\label{av-HSC-sol}
\delta_1(N)=\sqrt{\pi}\,\,\frac{\Gamma(N-1)}{\Gamma(N-\frac{3}{2})}
\end{equation}
diverging as $\sqrt{\pi N}$ when $N\to\infty$. \revv{This result was derived in \cite{DMS01}, together with the average degree of every vertex.} Similarly we deduce the recurrence relation 
\begin{equation*}
V_1(N+1)=\left[1+\frac{2}{2N-3}\right]V_1(N)+\frac{\delta_1(N)}{2N-3}-\left[\frac{\delta_1(N)}{2N-3}\right]^2
\end{equation*}
for the variance $\text{Var}[\delta_1]$ which we iterate starting with $V_1(3)=0$ to give
\begin{eqnarray}
\label{var-HSC-sol}
V_1(N) =4N-6-\delta_1(N)[1+\delta_1(N)]
\end{eqnarray}
diverging as 
\begin{equation*}
V_1(N)=(4-\pi) N-\sqrt{\pi N}+\frac{7\pi}{4}-6+O(N^{-1/2})
\end{equation*}
when $N\to\infty$. \revv{Equation~\eqref{var-HSC-sol} also follows from the exact distribution~\eqref{P-delta:sol} given below.} The growth laws for the average and the variance indicate that the degree $\delta_1$ is a non-self-averaging random quantity. 

The probability distribution $P(\delta_1,N)$ is expected to exhibit a scaling behavior in the $N\to\infty $ limit. More precisely, the growth laws $\delta_1(N)\sim \sqrt{N}$ and  $V_1(N)\sim N$ suggest that the probability distribution $P(\delta_1,N)$ approaches the scaling form
\begin{equation}
P(\delta_1,N) = \frac{1}{\sqrt{N}}\,\mathcal{P}(\Delta)
\end{equation}
in the scaling limit
\begin{equation}
\label{scaling}
N\to\infty, \quad \delta_1\to\infty, \quad \Delta=\frac{\delta_1}{\sqrt{N}} = \text{finite}
\end{equation}
\revv{The scaling function}
\begin{equation}
\revv{\mathcal{P}(\Delta)=\frac{\Delta}{2}\,e^{-\Delta^{2}/4}}
\label{eq:C_rayleigh}
\end{equation}
\revv{follows from the exact distribution~\eqref{P-delta:sol} below and agrees with $\langle\Delta\rangle=\sqrt{\pi}$ and $\langle\Delta^2\rangle-\langle\Delta\rangle^2=4-\pi$ obtained from~\eqref{av-HSC-sol}--\eqref{var-HSC-sol}.}

The probability distribution $P(\delta_1,N)$ satisfies the recurrence 
\begin{eqnarray}
\label{P-delta:HSC}
P(\delta_1,N+1) &=& \left(1-\frac{\delta_1}{2N-3}\right)P(\delta_1,N)\nonumber \\
&+&  \frac{\delta_1-1}{2N-3}\,P(\delta_1-1,N)
\end{eqnarray}
following from the stochastic evolution rule \eqref{delta-HSC-2}. The degree $\delta_1$ varies in the range
\begin{equation}
\label{delta:HSC}
2\leq \delta_1\leq N-1
\end{equation}

\revv{The recurrence~\eqref{P-delta:HSC} admits the exact solution~\cite{DMS01}}
\begin{equation}
\revv{P(\delta_1,N)=\frac{(\delta_1-1)\,(2N-\delta_1-4)!}{2^{\,N-1-\delta_1}\,(N-1-\delta_1)!\,(2N-5)!!}}
\label{P-delta:sol}
\end{equation}
\revv{for $2\le \delta_1\le N-1$. In particular, the extreme values occur with probabilities}
\begin{equation}
\revv{\begin{aligned}
P(2,N)&=\frac{1}{2N-5}\,,\\
P(N-1,N)&=\frac{\sqrt{\pi}}{2^{N-2}}\,\frac{\Gamma(N-1)}{\Gamma\!\left(N-\tfrac{3}{2}\right)}\,,
\end{aligned}}
\label{max:HSC}
\end{equation}
\revv{and, e.g., $P(3,N)=\frac{2}{2N-5}$, $P(4,N)=\frac{6(N-4)}{(2N-5)(2N-7)}$, $P(5,N)=\frac{8(N-5)}{(2N-5)(2N-7)}$.}

Let us probe the degree distribution. Relying again on the similarity in the growth of the degree of the primordial vertex and the maximal degree, $\delta_1\sim \delta_\text{max}$, and comparing \eqref{deg:max} with $\delta_1\sim\sqrt{N}$, we fix the exponent $\nu=3$ in \eqref{deg:tail}. \revv{This heuristic value agrees with the exact result~\cite{DMS01,Ginestra16-NGF,FIMS22}. The same argument applied to HSC($m$), where the link of a vertex grows as $N^{(m-1)/m}$, gives $\nu=2+1/(m-1)$, again the exact exponent~\cite{Ginestra16-NGF,FIMS22}. This supports the method used for the RRSC($m$) models in Sec.~\ref{sec:primordial}.} The minimal possible degree is $\delta_1=2$. As in Sec.~\ref{sec:degree}, we divide the set of all vertices into the set of $N_2$ vertices of degree two and the set of $N-N_2$ vertices of degree higher than two:
\begin{equation}
\label{vertices:HSC}
\underbrace{v_i \cdots}_{N_2} \quad \underbrace{v_I \cdots}_{N-N_2}
\end{equation}
We then divide the set of edges into two subsets:
\begin{equation}
\label{edges:HSC}
\underbrace{\{v_k,v_K\}\ldots}_{2N_2} \quad \underbrace{\{v_I,v_J\}\ldots}_{2N-3-2N_2}
\end{equation}

Adding a vertex, $N\to N+1$, leads to the evolution equation for the number of vertices of degree two
\begin{equation}
\label{N2:eq}
N_2\to 
\begin{cases}
N_2+1 &\text{prob} ~~1-\frac{2N_2}{2N-3}\\
N_2 &\text{prob} ~~\frac{2N_2}{2N-3}
\end{cases}
\end{equation}
Equation \eqref{N2:eq} is valid in the $N\geq 4$ range where the quantity $N_2$ can only increase. The (deterministic) initial condition reads 
\begin{equation}
\label{IC:2}
N_2 = 2  \qquad \text{when}\quad N=4
\end{equation}

Using \eqref{N2:eq} we deduce the recurrence 
\begin{equation}
\label{N2:av}
N_2(N+1)=\left[1-\frac{2}{2N-3}\right]N_2(N)+1
\end{equation}
for the average $N_2(N)=\langle N_2\rangle$. Starting with \eqref{IC:2} and iterating \eqref{N2:av} we obtain
\begin{equation}
\label{N2:av-sol}
N_2(N)=\frac{N^2-4N+6}{2N-5}
\end{equation}
Thus 50\% of vertices have the minimal degree two:
\begin{equation}
n_2 = \lim_{N\to\infty}\frac{N_2(N)}{N} = \frac{1}{2}\revv{,}
\end{equation}
\revv{in agreement with $P(2)=1/2$ from Refs.~\cite{DMS01,Ginestra16-NGF,FIMS22}; Eq.~\eqref{N2:av-sol} also follows from the exact finite-$N$ degree distribution of Ref.~\cite{DMS01}.}

\revv{To our knowledge, the fluctuations of $N_2$ have not been studied before.} We now show that $N_2$ is an asymptotically self-averaging random variable. Similarly to \eqref{N2:av} we deduce the recurrence for the variance $\langle\!\langle N_2^2\rangle\!\rangle$ which we solve to yield
\begin{equation}
\label{N2:var-sol}
\langle\!\langle N_2^2\rangle\!\rangle=\frac{(N-4)(2 N^3 - 12 N^2 + 25 N - 39)}{3(2N-7)(2N-5)^2}
\end{equation}
Thus $\langle\!\langle N_2^2\rangle\!\rangle \simeq  N/12$ and $N_2(N)\simeq N/2$ as $N\to\infty$, indicating that $N_2$ is an asymptotically self-averaging random variable since the relative magnitude of fluctuation $\sqrt{\langle\!\langle N_2^2\rangle\!\rangle}/N_2(N) \sim N^{-1/2}$ vanishes when $N\to\infty$. 

Using similar pedestrian calculations, one can compute the third cumulant
\begin{equation}
\label{N2-3:sol}
\langle\!\langle N_2^3\rangle\!\rangle=\frac{9(N-4)(7 N - 31)}{(2N-9)(2N-7)(2N-5)^3}
\end{equation}
An exact formula for the fourth cumulant is cumbersome, so we only mention the asymptotic $\langle\!\langle N_2^4\rangle\!\rangle \simeq -N/120$. The cumulants tend to grow linearly, $\langle\!\langle N_2^k\rangle\!\rangle \simeq C_k N$. For instance, $C_1=\frac{1}{2}, ~C_2=\frac{1}{12}, ~C_4=-\frac{1}{120}$. The third cumulant decays with $N$, see \eqref{N2-3:sol}, i.e.,  $C_3=0$.  To determine the amplitudes $C_k$, we use the stochastic evolution equation \eqref{N2:eq} and find
\begin{equation}
\label{N2-cum}
\left\langle e^{x N_2(N+1)}\right\rangle = e^x\left\langle e^{x N_2}\right\rangle+\frac{1-e^x}{N-3/2}\left\langle N_2 e^{x N_2}\right\rangle
\end{equation}
we write $N_2$ instead of $N_2(N)$; we use the complete notation only when the number of vertices is $N+1$.  It proves convenient to rewrite \eqref{N2-cum} as
\begin{equation}
\label{log:eq}
\left\langle e^{x N_2(N+1)}\right\rangle = e^x\left\langle e^{x N_2}\right\rangle+\frac{1-e^x}{N-3/2}\frac{d \left\langle e^{x N_2}\right\rangle}{d x}.
\end{equation}
The linear growth of the cumulants with $N$ suggests to define the (normalized) cumulant generating function via 
\begin{equation}
\label{C:def}
C(x) = \lim_{N\to\infty} \frac{\log \left\langle e^{x N_2}\right\rangle}{N-3/2}
\end{equation}
We thus recast \eqref{log:eq} into an ordinary differential equation (ODE)
\begin{equation}
\label{C:eq}
e^{C} = e^x +(1-e^x)\,\frac{d C}{d x}.
\end{equation}
The substitution $C = - \log u$ transforms \eqref{C:eq} into a linear (inhomogeneous) ODE which we solve subject to $u(0)=1$ corresponding to $C(0)=0$. Returning back to $C(x)$, we arrive at 
\begin{equation}
\label{C:sol}
C(x) = \log\,\frac{e^x-1}{x}
\end{equation}
The normalized cumulants, i.e., the amplitudes $C_k$, are encoded in the expansion 
\begin{equation}
\label{cum:exp}
C(x) = \sum_{k\geq 1} C_k\,\frac{x^k}{k!}
\end{equation}
Expanding $C(x)$ given by \eqref{C:sol} and comparing with \eqref{cum:exp} yields $C_k = k^{-1} B_k$ expressing the amplitudes via Bernoulli numbers $B_k$. (We use the convention $B_1=\frac{1}{2}$.) Thus $C_k=0$ for all odd $k\geq 3$. We also recover $C_k$ with $k\leq 4$ derived from pedestrian calculations. The next amplitudes are $C_6=\frac{1}{252}, ~C_8=-\frac{1}{240}, ~C_{10}=\frac{1}{132}$, etc.

The asymptotic behaviors of the cumulants suggest that the distribution $Q(N_2,N)$ of  the number of vertices of degree two is asymptotically Gaussian, characterized by the average and variance computed above:
\begin{equation}
\label{QN:Gauss}
Q(N_2,N)\simeq \sqrt{\frac{6}{\pi\,N}}\,\exp\!\left[-\frac{6(N_2-N/2)^2}{N}\right]
\end{equation}

To derive \eqref{QN:Gauss} and establish a few exact results for $Q(N_2,N)$  we first recast \eqref{N2:eq} into an exact recurrence 
\begin{eqnarray}
\label{QN2}
Q(N_2,N+1) &=& \left[1-\frac{2(N_2-1)}{2N-3}\right] Q(N_2-1,N)   \nonumber \\
&+& \frac{2N_2}{2N-3} Q(N_2,N) 
\end{eqnarray}
Iterating \eqref{QN2} and using the initial condition \eqref{IC:2}, that is, $Q(N_2,4)=\delta_{N_2,2}$, one can recurrrently calculate $Q(N_2,N)$ for $N>4$. We have not found an explicit general formula for $Q(N_2,N)$. Some exact results, e.g., explicit formulae for the extreme probabilities
\begin{subequations}
\begin{align}
\label{min:Q}
&Q(2,N)=2^{N-4}\frac{\Gamma(\frac{5}{2})}{\Gamma(N-\frac{3}{2})}\\
\label{max:Q}
& Q(N-2,N)=2^{4-N}\frac{\Gamma(\frac{5}{2})}{\Gamma(N-\frac{3}{2})}
\end{align}
\end{subequations}
readily follow from  \eqref{QN2}. 

To probe the asymptotic behavior of the probability distribution $Q(N_2, N)$, we employ a continuum description. Replacing the differences in \eqref{QN2} by derivatives and keeping only dominant terms, we recast the recurrence \eqref{QN2} into a PDE
\begin{equation}
\label{PDE}
\frac{\partial Q}{\partial N} = \frac{Q}{N} -\left[1-\frac{N_2}{N}\right]\frac{\partial Q}{\partial N_2}
+\frac{1}{2}\left[1-\frac{N_2}{N}\right]\frac{\partial^2 Q}{\partial N_2^2}
\end{equation}
It proves useful to transform the variables
\begin{equation}
\label{transform}
(N_2,N)\to (y=N_2-N/2, N)
\end{equation}
so that the peak of the distribution is at $y=0$. In the new variables \eqref{transform}, equation \eqref{PDE} becomes 
\begin{equation}
\label{PDE-y}
\frac{\partial Q}{\partial N} = \frac{Q}{N}+ \frac{y}{N}\,\frac{\partial Q}{\partial y}+\frac{1}{4}\,\frac{\partial^2 Q}{\partial y^2}
\end{equation}
where we kept again only the leading terms. The solution of \eqref{PDE-y} has a scaling form
\begin{equation}
\label{QNy}
Q(y,N) = \frac{1}{\sqrt{N}}\,\mathcal{Q}(z), \qquad  z =  \frac{y}{\sqrt{N}}
\end{equation}
By inserting \eqref{QNy} into \eqref{PDE-y} we obtain an ODE for the scaled distribution
\begin{equation}
\frac{d^2 \mathcal{Q}}{dz^2} + 6z\,\frac{d \mathcal{Q}}{dz}+6\mathcal{Q} =0
\end{equation}
The solution is the Gaussian distribution
\begin{equation}
\label{Gauss}
\mathcal{Q}(z) =\sqrt{\frac{3}{\pi}}\,  e^{-3 z^2}
\end{equation}
Recalling \eqref{transform} and \eqref{QNy}, we turn \eqref{Gauss} into \eqref{QN:Gauss}.

\bibliography{references-nets}

\end{document}